\documentclass[notitlepage,11pt]{article}
\usepackage{amssymb,amsmath,comment}
\catcode`\@=11
\@addtoreset{equation}{section}

\catcode`\@=12

\def\R{\mathbb{R}}
\def\C{\mathcal{C}}
\def\S{\mathbb{S}}
\def\f{\varphi}
\def\eps{\varepsilon}
\def\starstar{{2^{*\!*}}}
\def\uu{\underline u}
\def\div{{\rm div}}
\def\weak{{\,\rightharpoonup\,}}

\def\proof{\noindent{\textbf{Proof. }}}
\def\QED{\hfill {$\square$}\goodbreak \medskip}

\newtheorem{Theorem}{Theorem}[section]
\newtheorem{Lemma}[Theorem]{Lemma}
\newtheorem{Proposition}[Theorem]{Proposition}
\newtheorem{Corollary}[Theorem]{Corollary}
\newtheorem{Remark}[Theorem]{Remark}

\begin{document}
\title{\vspace{-10mm}
\Large{Caffarelli-Kohn-Nirenberg type inequalities for the weighted biharmonic operator: existence of extremal functions,\\{breaking} positivity and breaking symmetry}}


\author{Paolo Caldiroli \and Roberta Musina}



\date{}

\maketitle


\begin{abstract}
\footnotesize
We investigate Caffarelli-Kohn-Nirenberg type inequalities for the weighted biharmonic operator on
cones, both
under  Navier and Dirichlet boundary conditions.
Moreover, we study existence and qualitative properties of  extremal functions. In particular, we show that in some cases extremal functions do change sign; when the domain is the whole space, we prove some breaking symmetry phenomena.
\medskip

\medskip

\end{abstract}
%

\bigskip

\footnotesize

\noindent
{\bf Contents}

\noindent
\ref{S:Introduction}Introduction \dotfill \pageref{S:Introduction}


\noindent
\ref{S:Inequalities}\quad Inequalities \dotfill \pageref{S:Inequalities}


\ref{SS:proof}\quad Proof of Theorem \ref{T:CKN-Navier} \dotfill \pageref{SS:proof}

\ref{SS:Limiting}\quad Large and strict inequalities in the limiting case \dotfill \pageref{SS:Limiting}

\ref{SS:rad}\quad Radially symmetric functions \dotfill \pageref{SS:rad}

\ref{Ss:esti-Srad}\quad{Estimates on $S_{q}^{\mathrm{rad}}(\R^{n};\alpha)$} \dotfill 
\pageref{Ss:esti-Srad}

\noindent
\ref{S:existence}\quad Existence \dotfill \pageref{S:existence}

\ref{SS:eps-comp}\quad $\eps$-compactness \dotfill \pageref{SS:eps-comp}

\ref{SS:subcritical}\quad Proof of Theorem \ref{T:achieved1} \dotfill \pageref{SS:subcritical}

\ref{SS:critical}\quad Proof of Theorem \ref{T:achieved1critical} \dotfill \pageref{SS:critical}

\noindent
\ref{S:BP}\quad Breaking positivity \dotfill \pageref{S:BP}

\noindent
\ref{S:BS}\quad Breaking symmetry \dotfill \pageref{S:BS}

\noindent
\ref{S:Dbc}\quad Dirichlet boundary conditions \dotfill \pageref{S:Dbc}

\noindent
Appendix A: Remarks on a Brezis-Nirenberg type problem \dotfill \pageref{S:BN}

\noindent
Appendix B: Auxiliary results and open problems  \dotfill \pageref{S:Appendix}


\noindent
References \dotfill \pageref{References}

\normalsize

\section*{Introduction}\label{S:Introduction}

In this paper we study second order interpolation inequalities with weights being powers of the distance from the origin, and involving functions defined on dilation invariant domains. More precisely, for any regular domain $\Sigma$  in the unit sphere $\S^{n-1}$ we denote by $\mathcal{C}_{\Sigma}$ the cone
\begin{equation}
\label{eq:CSigma}
\mathcal{C}_{\Sigma}:  =  \left\{~r\sigma~|~r>0~,~\sigma\in\Sigma~\right\}.
\end{equation}
We are mainly interested in a class of inequalities of the form
\begin{equation}
\label{eq:intro-inequality}
\int_{\mathcal C_\Sigma}|x|^{\alpha}|\Delta u|^{2}dx\ge C\left(\int_{\mathcal C_\Sigma}
|x|^{-\beta}|u|^{q}dx\right)^{2/q}
\quad\textrm{for any $u\in C^{2}_{c}(\overline{\mathcal C_\Sigma}\setminus\{0\})$,}
\end{equation}
where $q>2$ and $\alpha\in\R$ are given parameters, and where
$C^{2}_{c}(\overline{\mathcal C_\Sigma}\setminus\{0\})$ is the space of functions in $C^2(\overline{\mathcal C_\Sigma})$ vanishing on $\partial\mathcal C_\Sigma$ and in a neighborhood of $0$ and of $\infty$. 
The best constant in (\ref{eq:intro-inequality}) is given by
\begin{equation}
\label{eq:bc}
S_q(\mathcal C_\Sigma;\alpha):  =  \inf_{\scriptstyle u\in C^{2}_{c}(
\overline{\mathcal C_\Sigma}\setminus\{0\})
\atop\scriptstyle u\ne 0}\frac{\displaystyle
\int_{\mathcal C_\Sigma}|x|^{\alpha}|\Delta u|^{2}dx}
{\displaystyle\left(\int_{\mathcal C_\Sigma}|x|^{-\beta}|u|^{q}dx\right)^{2/q}}~\!.
\end{equation}
A simple rescaling argument shows that $S_q(\mathcal C_\Sigma;\alpha)$ vanishes unless
\begin{equation}
\label{eq:beta}
\beta  =  n-q~\!\frac{n-4+\alpha}{2}~\!.
\end{equation}
Therefore, from now on we will assume that (\ref{eq:beta}) holds.

If $n\ge 5$ then the Sobolev embedding theorem implies that a necessary condition to have
$S_q(\mathcal C_\Sigma;\alpha)>0$ is that $q\le 2^{*\!*}$, where $2^{*\!*}$ is the critical Sobolev exponent:
$$
2^{*\!*}  =  \frac{2n}{n-4}~\!.
$$
Our goal is to  estimate the best constant $S_q(\mathcal C_\Sigma;\alpha)$
under the above assumptions on $\beta$ and $q$. 
Moreover, we will study existence and qualitative properties of functions achieving $S_q(\mathcal C_\Sigma;\alpha)$
on a suitable function space.

Let us notice that (\ref{eq:intro-inequality}) can not be obtained by iterating
the first order Caffarelli-Kohn-Nirenberg inequalities in \cite{CKN}, see
Remark  \ref{R:Lin}. We quote \cite{Lin86}, \cite{SzWa} and references there-in, for a related interpolation inequality due to 
C.S. Lin.

If $q  =  2$ then (\ref{eq:beta}) gives $\beta  =  \alpha-4$ and (\ref{eq:intro-inequality}) becomes
\begin{equation}
\label{eq:Rellich}
\int_{\mathcal C_\Sigma}|x|^{\alpha}|\Delta u|^{2}dx\ge C\int_{\mathcal C_\Sigma}
|x|^{\alpha-4}|u|^{2}dx
\quad\textrm{for any $u\in C^{2}_{c}(\overline{\mathcal C_\Sigma}\setminus\{0\})$.}
\end{equation}
A first version of this inequality has been introduced by F. Rellich in 1953 (see \cite{Rel54} and \cite{Rel69}) in case $\mathcal C_\Sigma  =  \R^{n}\setminus\{0\}$ and $\alpha  =  0$. 
For general cones $\mathcal C_\Sigma$ and parameters $\alpha\in\R$
we refer to \cite{CM1}, where it is proved that the best constant in (\ref{eq:Rellich})  is 
exactly the square of the distance of $-\gamma_\alpha$ from the Dirichlet spectrum $\Lambda(\Sigma)$ of the Laplace-Beltrami operator on $\Sigma$, where
\begin{equation}
\label{eq:gamma}
\gamma_{\alpha}  =  \left(\frac{n-2}{2}\right)^{2}-\left(\frac{\alpha-2}{2}\right)^{2}.
\end{equation}
For instance, taking $\Sigma  =  \S^{n-1}$ or $\Sigma = $ half-sphere we have
\begin{eqnarray*}
S_2(\R^{n}\setminus\{0\};\alpha)
 & = & \min_{k\in\mathbb{N}\cup\{0\}}\left|\gamma_{\alpha}+k(n-2+k)\right|^{2}\\
S_2(\R^{n}_+;\alpha)
 & = & \min_{k\in\mathbb{N}}\left|\gamma_{\alpha}+k(n-2+k)\right|^{2}~\!,
\end{eqnarray*}
where $\R^n_+$ denotes any homogeneous half-space.

In our first theorem we show that $S_q(\C_\Sigma,\alpha)>0$
whenever the best constant in the weighted Rellich inequality is positive.

\begin{Theorem}
\label{T:CKN-Navier}
Let $\alpha\in\R$ and let $\Sigma\subseteq\S^{n-1}$ be a domain of class $C^2$. Let $q>2$ be a given exponent, and assume that $q\le 2^{*\!*}$ if $n\ge 5$. Then $S_q(\mathcal C_\Sigma;\alpha)>0$ if and only if $-\gamma_{\alpha}\notin\Lambda(\Sigma)$.
\end{Theorem}

If $-\gamma_{\alpha}$ is not a Dirichlet eigenvalue on $\Sigma$ then we can define the Hilbert space $\mathcal{N}^{2}(\mathcal{C}_{\Sigma};\alpha)$ as the completion of $C^{2}_{c}(\overline{\mathcal{C}_{\Sigma}}\setminus\{0\})$ with respect to the norm
\begin{equation}
\label{eq:norm}
\|u\|_{2,\alpha}   =  \left(\int_{\mathcal{C}_{\Sigma}}|x|^{\alpha}|\Delta u|^{2}dx\right)^{1/2}~\!.
\end{equation}
If $n\ge 5$ and $\alpha=0$ then $\mathcal{N}^{2}(\R^n\setminus\{0\};0) =\mathcal{D}^{2}(\R^{n})$,
see Remark \ref{R:density}. In general, it holds that
\begin{equation}
\label{eq:minimization1}
S_q(\mathcal C_\Sigma;\alpha)  =  \inf_{\scriptstyle u\in \mathcal{N}^{2}(\mathcal{C}_{\Sigma};\alpha)\atop\scriptstyle u\ne 0}\frac{\displaystyle\int_{\mathcal C_\Sigma}|x|^{\alpha}|\Delta u|^{2}dx}
{\displaystyle\left(\int_{\mathcal C_\Sigma}|x|^{-\beta}|u|^{q}dx\right)^{2/q}}~\!.
\end{equation}  

In the rest of the paper we study the existence of extremals for $S_q(\mathcal C_\Sigma;\alpha)$
and their qualitative properties.

When $q  =  2$ it was shown in \cite{CM1} that the best constant $S_2(\C_\Sigma;\alpha)$ is never attained in $\mathcal{N}^{2}(\mathcal{C}_{\Sigma};\alpha)$. Another
remarkable case is
  $\Sigma  =  \S^{n-1}$, 
$n\ge 5$ and $q  =  2^{*\!*}$. Then $\mathcal C_\Sigma  =  \R^{n}\setminus\{0\}$,
$\beta  =  0$ and
$S_{2^{*\!*}}(\R^{n}\setminus\{0\};0)$ equals the Sobolev constant 
\begin{equation}
\label{eq:Sob}
S^{*\!*}  =  \inf_{\scriptstyle u\in\mathcal D^{2}(\R^n)
\atop\scriptstyle u\ne 0}\frac{\displaystyle
\int_{\R^{n}}|\Delta u|^{2}dx}
{\displaystyle\left(\int_{\R^{n}}|u|^{2^{*\!*}}dx\right)^{2/2^{*\!*}}}~.
\end{equation}
It is well known that the best constant $S^{*\!*}$ is achieved by 
an explicitly known radially symmetric and positive function, see for instance \cite{Sw}.

In the next results we study the attainability of
$S_q(\mathcal C_\Sigma;\alpha)$. 
By standard arguments, extremals for $S_q(\mathcal C_\Sigma;\alpha)$ are, up to a Lagrange multiplier, ground state solutions of equation
\begin{equation}
\label{eq:eq}
\Delta(|x|^{\alpha}\Delta u)  =  |x|^{-\beta}|u|^{q-2}u\quad\textrm{in }\mathcal{C}_{\Sigma}
\end{equation}
under Navier boundary conditions $u  =  \Delta u  =  0$ on $\partial\mathcal{C}_{\Sigma}$, in case $\Sigma$ is properly contained in $\S^{n-1}$.

Notice that the minimization problem
(\ref{eq:Sob}) is noncompact, due to the action of the group of dilations in $\R^n$.
However, when $q$ is subcritical the infimum $S_q(\mathcal C_\Sigma;\alpha)$ is always
achieved.

\begin{Theorem}
\label{T:achieved1}
Let $q>2$ be a given exponent such that $q< 2^{*\!*}$ if $n\ge 5$.  
Let $\Sigma$ be a domain in $\S^{n-1}$ of class $C^2$, with $-\gamma_{\alpha}\notin\Lambda(\Sigma)$. 
Then $S_q(\mathcal C_\Sigma;\alpha)$ is achieved in 
$\mathcal{N}^{2}(\mathcal{C}_{\Sigma};\alpha)$.
\end{Theorem}

When $n\ge 5$ and $q  =  2^{*\!*}$ it holds that $S_{2^{*\!*}}(\mathcal C_\Sigma;\alpha)\le S^{*\!*}$ for any cone $\mathcal C_\Sigma$ and for any admissible exponent $\alpha$, see Proposition \ref{P:large}.
In this case the group of translations in $\R^n$
may produce lack of compactness and nonexistence phenomena.
As usual, the strict inequality guarantees the compactness of all minimizing sequences.

\begin{Theorem}
\label{T:achieved1critical}
Let $n\ge 5$ and let $\Sigma$ be a domain in $\S^{n-1}$ of class $C^2$. Assume that $-\gamma_{\alpha}\notin\Lambda(\Sigma)$. If $S_{2^{*\!*}}(\mathcal C_\Sigma;\alpha)<S^{*\!*}$ then $S_{2^{*\!*}}(\mathcal C_\Sigma;\alpha)$ is achieved in $\mathcal{N}^{2}(\mathcal{C}_{\Sigma};\alpha)$.
\end{Theorem}

The above stated theorems constitute the second order version of well known results related to the classical Caffarelli-Kohn-Nirenberg inequalities \cite{CKN} for first order operators, see for instance
\cite{CalMus}, \cite{CatWan} and the references therein.

However, when we push further the study of minimization problems (\ref{eq:minimization1}), some meaningful differences appear. Firstly we can show that in the case of critical exponent the strict inequality $S_{2^{*\!*}}(\mathcal C_\Sigma;\alpha)<S^{*\!*}$ holds in the following cases.

\begin{Theorem}
If $n\ge 6$ and $|\alpha-2|>2$ then $S_{2^{*\!*}}(\mathcal{C}_{\Sigma};\alpha)<S^{*\!*}$ for every $\Sigma\subseteq\S^{n-1}$.\\
If $n  =  5$ and $2<|\alpha-2|<\sqrt{13}$ then $S_{2^{*\!*}}(\R^{5}\setminus\{0\};\alpha)<S^{*\!*}$.
\end{Theorem}

The previous result is discussed separately for dimensions $n\ge 6$ in Theorem \ref{T:minore}, whereas for $n  =  5$ is a special case of an estimate proved in Theorem \ref{T:critico-RN}. 

The difference between the case $n  =  5$ and $n\ge 6$ seems to be not purely techical. There is indeed a deep connection between the  validity of the strict inequality
$S_{2^{*\!*}}(\mathcal{C}_{\Sigma};\alpha)<S^{*\!*}$ and the existence of ground state
solutions for  the following
Dirichlet problem:
\begin{equation}
\label{eq:BN-problem}
\begin{cases}
\Delta^2 u+\lambda\Delta u  =  |u|^{\starstar-2}u&\textrm{in $B$}\\
u  =  |\nabla u|  =  0&\textrm{on $\partial B$.}
\end{cases}
\end{equation}
Here  $B\subset \R^n$ is the unit ball and $\lambda$ is a given real parameter.
As a by-product of our computations we can prove a
Brezis-Nirenberg type result  for problem (\ref{eq:BN-problem}) in the spirit of the celebrated
paper \cite{BreNir}, see 
Appendix \ref{S:BN}. By adapting a terminology which has been introduced
by Pucci and Serrin in \cite{PuSe}, we can assert that $n  =  5$ is the unique {\bf weakly critical dimension} for problem (\ref{eq:BN-problem}).

When $\mathcal{C}_{\Sigma}  =  \R^{n}\setminus\{0\}$ \textbf{breaking symmetry} can be observed as well. In particular, from the results in Section \ref{S:BS} it follows that minimizers for $S_q(\mathcal C_\Sigma;\alpha)$ may be not radially symmetric. In Theorems \ref{T:break-symm1} and \ref{T:bs} we show that breaking symmetry occurs, for instance, when $-\gamma_\alpha$ is close  to a Dirichlet eigenvalue on the sphere or when $|\alpha|$ is large enough.

Even in correspondence of the critical exponent, breaking symmetry occurs: for $|\alpha|$ large enough there exist minimizers both for $S_{2^{*\!*}}(\R^{n}\setminus\{0\};\alpha)$ and for the corresponding radial best constant $S_{2^{*\!*}}^{\mathrm{rad}}(\R^{n};\alpha)$, defined in (\ref{eq:minimization1radial}), and the minimizers are different as $S_{2^{*\!*}}(\R^{n}\setminus\{0\};\alpha)<S_{2^{*\!*}}^{\mathrm{rad}}(\R^{n};\alpha)$. A similar breaking symmetry phenomenon
does not occur, for instance, in dealing with critical exponents in first-order
Caffarelli-Kohn-Nirenberg inequalities: in that case, the best constant is not achieved, or all the
minimizers are radially symmetric. We refer to \cite{CatWan}, \cite{FelSch},
\cite{DELT} for breaking symmetry in first order
Caffarelli-Kohn-Nirenberg inequalities.

Another striking difference with respect to similar first order problems, is a \textbf{breaking positivity} phenomenon. Indeed, in Section \ref{S:BP} we show that, in general, no extremal for (\ref{eq:minimization1}) has constant sign,
see Theorem \ref{T:break-pos}.

In Section \ref{S:Dbc} we take $\Sigma$ to be a proper domain in the sphere and 
we deal with the infimum
\begin{equation}\label{eq:dirichlet}
S_q^D(\mathcal C_\Sigma;\alpha):  =  \inf_{\scriptstyle u\in C^{2}_{c}(\mathcal C_\Sigma)\atop\scriptstyle u\ne 0}\frac{\displaystyle\int_{\mathcal C_\Sigma}|x|^{\alpha}|\Delta u|^{2}dx}{\displaystyle\left(\int_{\mathcal C_\Sigma}|x|^{-\beta}|u|^{q}dx\right)^{2/q}}~\!.
\end{equation}
Differently from the Navier case, it turns out that $S_q^D(\mathcal C_\Sigma;\alpha)$ is always positive, whenever $\Sigma$ has compact closure in $\S^{n-1}$, with no restriction on $\alpha$. We also show existence of extremals (see Theorem \ref{T:Dir}) which give rise to solutions of (\ref{eq:eq}) satisfying Dirichlet boundary conditions $u=|\nabla u|=0$ on $\partial\mathcal{C}_{\Sigma}$.


\section{Inequalities}\label{S:Inequalities}

In this Section we prove Theorem \ref{T:CKN-Navier} and some related results.
We start by noticing that 
\begin{equation}
\label{eq:alpha4alpha}
S_q(\C_\Sigma;\alpha)  =  S_q(\C_\Sigma;4-\alpha)~\!
\quad\textrm{and}\quad
S^D_q(\C_\Sigma;\alpha)  =  S^D_q(\C_\Sigma;4-\alpha)~\!.
\end{equation}
To check (\ref{eq:alpha4alpha}) use (as in \cite{CM1}, where $q=2$ is assumed) the
transform $u\mapsto \hat{u}$ given by
$$
\hat{u}(x)  =  |x|^{2-n}u(|x|^{-2}x).
$$ 
The proof of Theorem \ref{T:CKN-Navier} is based on the Emden-Fowler transform
$u\mapsto w=Tu$, defined by
\begin{equation}
\label{eq:EF}
u(x)  =  |x|^{\frac{4-n-\alpha}{2}}~\!w\left(-\log|x|,\frac{x}{|x|}\right).
\end{equation}
Notice that $T$ maps functions $u\colon\overline{\mathcal{C}_{\Sigma}}\setminus\{0\}\to\R$ into
functions $w$ on the cylinder
$$
\mathcal{Z}_{\Sigma}:  =  \{(s,\sigma)\in\R\times\S^{n-1}~|~s\in\R,~\sigma\in\Sigma\}.
$$
%
%
In \cite{CM1} it is noticed that for 
every $u\in C^{2}_{c}(\overline{\mathcal C_\Sigma}\setminus\{0\})$ one has
$w\in C^2_c(\overline{\mathcal Z_\Sigma})$ and
\begin{gather}
\label{eq:Lp}
\int_{\mathcal{C}_{\Sigma}}|x|^{-\beta}|u|^q dx  =   
\int_{\mathcal Z_\Sigma}|w|^q dsd\sigma\\
\label{eq:Delta2}
\int_{\mathcal{C}_{\Sigma}}|x|^{\alpha}|\Delta u|^{2}dx  =  \int_{\mathcal{Z}_{\Sigma}}\!\left(|L_{\alpha}w|^{2}+|w_{ss}|^{2}+2|\nabla_{\sigma}w_{s}|^{2}+2\overline\gamma_{\alpha}|w_{s}|^{2}\right)dsd\sigma~\!,
\end{gather}
where 
\begin{equation}
\label{eq:gbar}
L_{\alpha}w  =  -\Delta_{\sigma}w+\gamma_{\alpha}w~,\quad
\overline\gamma_{\alpha}  =  \left(\frac{n-2}{2}\right)^{2}+\left(\frac{\alpha-2}{2}\right)^{2}~\!,
\end{equation}
and $\gamma_{\alpha}$ is defined in (\ref{eq:gamma}).
For every $\gamma\in\R$ we introduce also the value
\begin{equation}
\label{eq:mN}
m_{N}(\Sigma;\gamma)  =  \inf_{\scriptstyle\varphi\in H^{2}\cap H^{1}_{0}(\Sigma)\atop\scriptstyle\varphi\ne 0}\frac{\displaystyle\int_{\Sigma}|-\Delta_{\sigma}\varphi+\gamma\varphi|^{2}d\sigma}{\displaystyle\int_{\Sigma}\varphi^{2}d\sigma}~\!.
\end{equation}
The following facts hold (see Proposition 1.1 and Theorem 2.1 in \cite{CM1}).
\begin{Lemma}\label{L:aux}
\begin{itemize}
\item[$(i)$] 
For every $\gamma\in\R$ one has that $m_{N}(\Sigma;\gamma)  =  \mathrm{dist}(-\gamma,\Lambda(\Sigma))^{2}$. Moreover $\varphi\in H^{2}\cap H^{1}_{0}(\Sigma)$ is a minimizer for $m_{N}(\Sigma;\gamma)$ if and only if $\varphi$ is an eigenfunction of $-\Delta_{\sigma}$ relative to the eigenvalue achieving the minimal distance of $-\gamma$ from $\Lambda(\Sigma)$.
\item[$(ii)$] 
For every $\alpha\in\R$ one has that $S_{2}(\mathcal{C}_{\Sigma};\alpha)  =  m_{N}(\Sigma;\gamma_{\alpha})$ with $\gamma_{\alpha}$ given by (\ref{eq:gamma}).
\end{itemize}
\end{Lemma}


\subsection{Proof of Theorem \ref{T:CKN-Navier}}\label{SS:proof}
Assume that $-\gamma_\alpha\in \Lambda(\Sigma)$, and take a nontrivial $\f\in H^1_0(\Sigma)$ in the kernel of the operator $L_{\alpha}$. Test $S_{q}(\mathcal C_\Sigma;\alpha)$ with 
$$
u(x)  =  |x|^{\frac{4-n-\alpha}{2}}~\!\eta(-\log|x|)~\!\f\left(\frac{x}{|x|}\right)~\!,
$$
where $\eta\in C^\infty_c(\R)$, $\eta\neq 0$  is an arbitrary function. Using (\ref{eq:Lp})--(\ref{eq:Delta2}) we readily get
$$
S_{q}(\mathcal C_\Sigma;\alpha)\le C_\f~\!\frac{\displaystyle\int_{-\infty}^\infty|\eta''|^2ds+\int_{-\infty}^\infty|\eta'|^2 ds}
{\left(\displaystyle\int_{-\infty}^\infty|\eta|^q ds\right)^{2/q}}~\!,
$$
where the constant $C_\f>0$ does not depend on $\eta$. Thus $S_{q}(\mathcal C_\Sigma;\alpha)  =  0$, by a simple rescaling argument.

Next, assume that $-\gamma_\alpha\notin \Lambda(\Sigma)$. By the results in \cite{CM1}, it turns out that the space $H^{2}\cap H^{1}_{0}(\mathcal Z_\Sigma)$ has an equivalent norm given by
$$
\|w\|_{H^{2}\cap H^{1}_{0}(\mathcal Z_\Sigma;\alpha)}^{2}  =  \int_{\mathcal{Z}_{\Sigma}}\!\left(|L_{\alpha}w|^{2}+|w_{ss}|^{2}+2|\nabla_{\sigma}w_{s}|^{2}+2\overline\gamma_{\alpha}|w_{s}|^{2}\right)dsd\sigma.
$$
Moreover, the operator $T$ is an isomorphism between the spaces $\mathcal{N}^{2}(\mathcal{C}_{\Sigma};\alpha)$ and $H^{2}\cap H^{1}_{0}(\mathcal Z_\Sigma)$ and (\ref{eq:Lp})--(\ref{eq:Delta2}) hold for every $u\in \mathcal{N}^{2}(\mathcal{C}_{\Sigma};\alpha)$.
In addition, thanks to the Sobolev embedding theorem for $H^{2}(\mathcal Z_\Sigma)$ (see \cite{Ad}) and by (\ref{eq:Lp})--(\ref{eq:Delta2}), we infer that $\mathcal{N}^{2}(\mathcal{C}_{\Sigma};\alpha)$ is continuously embedded into $L^{q}(\mathcal{C}_{\Sigma};|x|^{-\beta}dx)$, namely, $S_{q}(\mathcal C_\Sigma;\alpha)>0$. 
\QED

\begin{Remark}
\label{R:density}
Let $\alpha\in\R$, $q>2$ with $q\le 2^{*\!*}$ if $n\ge 5$, and 
$\beta  =  n-q~\!\frac{n-4+\alpha}{2}$. If $n>4-\alpha$ then $C^{2}_{c}(\R^{n})\subset L^{q}(\R^{n};|x|^{-\beta}dx)$ and
\begin{equation}
\label{eq:SqRnalpha}
S_q(\R^{n}\setminus\{0\};\alpha)  =  \inf_{\scriptstyle u\in C^2_c(\R^{n})\atop\scriptstyle u\ne 0}\frac{\displaystyle
\int_{\R^{n}}|x|^{\alpha}|\Delta u|^{2}dx}
{\displaystyle\left(\int_{\R^{n}}|x|^{-\beta}|u|^{q}dx\right)^{2/q}}~\!.
\end{equation}
Moreover if $-\gamma_{\alpha}\not\in\Lambda(\S^{n-1})$ and $n>4-\alpha$ then $C^2_c(\R^{n})$ is dense in $\mathcal{N}^{2}(\R^{n}\setminus\{0\};\alpha)$. 
These facts can be proved in a standard way. 
\end{Remark}

\begin{Remark}
\label{R:Lin}
C.S. Lin proved in \cite{Lin86} several interpolation inequalities involving weighted
$L^p$ norms of the derivatives of functions $u\in C^\infty_c(\R^n)$. In particular, in case
$n\ge 5$ he proved that
for {any} $\alpha\in\R$ and for any $q\in(2,\starstar]$, there exists $C_L>0$ such that
$$
\max_{i,j=1,...,n}\int_{\R^n}|x|^{\alpha}|\partial_{ij}u|^{2}dx\ge C_L\left(\int_{\R^n}
|x|^{-\beta}|u|^{q}dx\right)^{2/q}
\quad\textrm{for any $u\in C^{2}_{c}({\R^n\setminus\{0\}})$,}
$$
where $\beta$ is given by (\ref{eq:beta}). Clearly, the best constant  $C_L$ controls
$S_q(\R^n\setminus\{0\};\alpha)$ from above,
but it always
happen that
$$
0=S_q(\R^n\setminus\{0\};\alpha) <C_L(\alpha)
$$
when $-\gamma_\alpha=k(n-2+k)$ for some positive integer $k$. More precisely,
the functions
$$
u\mapsto\max_{i,j=1,...,n}~\int_{\R^{n}}|x|^\alpha|\partial_{ij}u|^{2}dx~,
\quad
u\mapsto \int_{\R^{n}}|x|^\alpha|\Delta u|^2~\!dx
$$
define two equivalent norms in $C^\infty_c(\R^n)$ {if and only if} $-\gamma_\alpha\notin
\Lambda(\S^{n-1})$. Notice that the present remark improves Lemma 3.1
in \cite{SzWa} (in case $p=k=2$), where $4-n<\alpha\le 0$ is assumed.
\end{Remark}

\subsection{Large and strict inequalities in the limiting case}\label{SS:Limiting}
In this subsection we take $n\ge 5$ and $q  =  \starstar$. Let
$S_\starstar(\mathcal C_\Sigma;\alpha)$, $S^D_\starstar(\mathcal C_\Sigma;\alpha)$ be the infima defined in
(\ref{eq:bc}),
(\ref{eq:dirichlet}) respectively. In particular,
$$
S^D_\starstar(\mathcal C_\Sigma;\alpha)  =  \inf_{\scriptstyle u\in C^{2}_{c}(
\mathcal C_\Sigma)
\atop\scriptstyle u\ne 0}\frac{\displaystyle
\int_{\mathcal C_\Sigma}|x|^{\alpha}|\Delta u|^{2}dx}
{\displaystyle\left(\int_{\mathcal C_\Sigma}|x|^{\frac{n\alpha}{n-4}}|u|^{\starstar}dx\right)^{2/\starstar}}~.
$$

\begin{Proposition}
\label{P:large}
Let $\Sigma$ be a domain in $\S^{n-1}$ of class $C^2$, $n\ge 5$, and let $\alpha\in\R$. Then
$$
S_\starstar(\mathcal C_\Sigma;\alpha)\le S^D_\starstar(\mathcal C_\Sigma;\alpha)\le S^{*\!*}~\!,
$$
where $S^{*\!*}$ is the Sobolev constant, given by (\ref{eq:Sob}).
\end{Proposition}

\proof
The first inequality is trivial. To prove that $S^D_\starstar(\mathcal C_\Sigma;\alpha)\le S^{*\!*}$ we
fix a point $x_0\in \C_\Sigma$. For an arbitrary $u\in C^2_c(\R^n)$, $u\neq 0$ 
and for any integer $h>0$ we put
$$
u_h(x)  =  u(h(x-x_0)).
$$
If $h$ is large enough then the support of $u_h$ is compactly contained in $\C_\Sigma$, and hence
\begin{eqnarray*}
S^D_\starstar(\mathcal C_\Sigma;\alpha)&\le&
\frac{\displaystyle
\int_{\mathcal C_\Sigma}|x|^{\alpha}|\Delta u_h|^{2}dx}
{\displaystyle\left(\int_{\mathcal C_\Sigma}|x|^{\frac{n\alpha}{n-4}}|u_h|^{\starstar}dx\right)^{2/\starstar}}
=
\frac{\displaystyle
\int_{\R^n}\left|\frac{y}{h}+x_0\right|^{\alpha}|\Delta u|^{2}dx}
{\displaystyle\left(\int_{\R^n}\left|\frac{y}{h}+x_0\right|^{\frac{n\alpha}{n-4}}|u|^{\starstar}dx\right)
^{2/\starstar}}
\\
&  =  &
\frac{\displaystyle
\int_{\R^n}|\Delta u|^{2}dx}
{\displaystyle\left(\int_{\R^n}|u|^{\starstar}dx\right)^{2/\starstar}}+o(1)
\end{eqnarray*}
as $h\to \infty$. Since $u$ was arbitrarily chosen, the conclusion follows.
\QED

As concerns the validity of the strict inequality $S_{2^{*\!*}}^{D}(\C_\Sigma;\alpha)<S^{*\!*}$ we have the following result.

\begin{Theorem}\label{T:minore}
If $n\ge 6$ and $|\alpha-2|>2$ then $S_{2^{*\!*}}^{D}(\mathcal{C}_{\Sigma};\alpha)<S^{*\!*}$ for every $\Sigma\subset\S^{n-1}$.
\end{Theorem}

\proof
Let $a  =  -{\alpha}/{2}$. We notice that $|\alpha-2|>2$ is equivalent to say that $C_{a}:={a(a+2)(n-2)}/n>0$. 
By Lemma \ref{L:radial-esti} in Appendix \ref{S:Appendix}, there exists $T_{a}\in(0,1)$ such that 
for $0<t\le T_a$ and for every radial mapping $u\in C^{2}_{c}(B)$, where $B$ is the unit open ball in $\R^{n}$, one has
\begin{equation}
\label{eq:vau3}
\int|tx+e|^{-2a}\left|\Delta\left(|tx+e|^{a}u\right)\right|^{2}\le 
\int|\Delta u|^{2}-C_{a}t^{2}\int|\nabla u|^{2}.
\end{equation}
Fix a point $e\in\Sigma$. Let $t_{0}>0$ be such that $e+t_{0}B\subset\mathcal{C}_{\Sigma}$
and put $t=\frac{1}{2}~\!\min\{t_{0},T_a\}$.
By Lemma \ref{L:BN-esti} in Appendix \ref{S:BN}, there exists a radially 
symmetric function $u\in C^{2}_{c}(B)$ such that 
\begin{equation}
\label{eq:BN-esti}
\int|\Delta u|^{2}-C_{a}t^{2}\int|\nabla u|^{2}<S^{*\!*}\left(\int |u|^\starstar\right)^{2/2^{*\!*}}.
\end{equation}
Define
$$
v(x)  =  |x|^{-\frac{\alpha}{2}}u\left(\frac{x-e}{t}\right)
$$
and notice that $v\in C^{2}_{c}(\mathcal{C}_{\Sigma})$ verifies
\begin{equation*}
\int|x|^{\frac{n\alpha}{n-4}}|v|^{2^{*\!*}}  =  t^{n}\int|u|^{2^{*\!*}}~~\textrm{and}~~
\int|x|^{\alpha}|\Delta v|^{2}  =  t^{n-4}\int|tx+e|^{-2a}\left|\Delta\left(|tx+e|^{a}u\right)\right|^{2}.
\end{equation*}
Then, by (\ref{eq:vau3}) and (\ref{eq:BN-esti})
\begin{equation*}
\begin{split}
S_{2^{*\!*}}(\mathcal{C}_{\Sigma};\alpha)&\le\frac{\displaystyle\int|x|^{\alpha}|\Delta v|^{2}}{\displaystyle\left(\int|x|^{\frac{n\alpha}{n-4}}|v|^{2^{*\!*}}\right)^{2/2^{*\!*}}}  =  \frac{\displaystyle\int|tx+e|^{-2a}\left|\Delta\left(|tx+e|^{a}u\right)\right|^{2}}{\displaystyle\left(\int|u|^{2^{*\!*}}\right)^{2/2^{*\!*}}}\\
&\le\frac{\displaystyle\int|\Delta u|^{2}-C_{a}t^{2}\int|\nabla u|^{2}}{\displaystyle\left(\int|u|^{2^{*\!*}}\right)^{2/2^{*\!*}}}<S^{*\!*}.
\end{split}
\end{equation*}
\vspace{-25pt}\par
\hfill\QED

In dimension $n  =  5$ we have a partial result on the whole space.

\begin{Theorem}\label{T:critico-R5}
If $2<|\alpha-2|<\sqrt{13}$ then $S_{2^{*\!*}}(\R^{5}\setminus\{0\};\alpha)<S^{*\!*}$.
\end{Theorem}

Theorem \ref{T:critico-R5} is a special case of an estimate which will be proved in the next subsection (see Theorem \ref{T:critico-RN}).

\subsection{Radially symmetric functions}\label{SS:rad}

For any $\alpha\in \R$ and $q\ge 2$ we define
 \begin{equation}
 \label{eq:minimization1radial}
S_q^{\rm rad}(\R^{n};\alpha):  =  \inf_{\scriptstyle u\in C^2_c({\R^{n}}\setminus\{0\})\atop\scriptstyle u\ne 0,~u  =  u(|x|)}
\frac{\displaystyle
\int_{\R^{n}}|x|^{\alpha}|\Delta u|^{2}dx}
{\displaystyle\left(\int_{\R^{n}}|x|^{-\beta}|u|^{q}dx\right)^{2/q}}\end{equation}
where $\beta  =  n-q\frac{n-4+{\alpha}}{2}$.
Notice that there is no upper bound on $q$ even in large dimensions.
Arguing as for (\ref{eq:alpha4alpha}), one can easily check that
$$
S_q^{\rm rad}(\R^n;\alpha)  =  S_q^{\rm rad}(\R^n;4-\alpha)~\!.
$$
In case $q  =  2$, it was proved in \cite{CM1} that
\begin{equation}
\label{eq:inf-radial}
S_2^{\rm rad}(\R^n;\alpha)  = \gamma_\alpha^2
  =~\!  \frac{(n-4+\alpha)^2(n-\alpha)^2}{16}~\!.
\end{equation}
In particular, if $\alpha\neq 4-n$ and $\alpha\neq n$ then $S_2^{\rm rad}(\R^n;\alpha)> 0$,
and we can suitably define a Hilbert space of radially symmetric functions 
$\mathcal{N}^{2}_{\rm rad}(\R^n;\alpha)$ endowed with the norm (\ref{eq:norm}).

The next theorem provides a second order Caffarelli-Kohn-Nirenberg type inequality 
for radially symmetric maps. We only need to assume $\gamma_\alpha\neq 0$. In particular,
 $q$ can be supercritical and 
$-\gamma_{\alpha}$ might be
a Dirichlet eigenvalue on the sphere.

\begin{Theorem}
\label{T:CKNradial}
Let $q>2$ be a given exponent. Then $S_q^{\rm rad}(\R^{n};\alpha)>0$
if and only if $\alpha\not\in\{4-n,n\}$. 
\end{Theorem}

\proof
To any radial function $u\in C^2_c(\R^{n}\setminus\{0\})$ we associate a function $w\in C^2_c(\R)$ via the Emden-Fowler transform defined in (\ref{eq:EF}).
Thus
\begin{gather}
\label{eq:uwrad}
u(x)  =  |x|^{\frac{4-n-\alpha}{2}}w(-\log|x|)\\
\nonumber
\int_{\R^{n}}|x|^{-\beta}|u|^q dx  =  \omega_{n}\int_{-\infty}^\infty|w|^q ds\\
\nonumber
\displaystyle\int_{\R^{n}}|x|^{\alpha}|\Delta u|^{2}dx  =  
\omega_{n} \displaystyle\int_{-\infty}^\infty
\left(|w''|^{2}+2\overline{\gamma}_{\alpha}|w'|^{2}+\gamma_{\alpha}^{2}|w|^{2}\right)ds
\end{gather}
where $\omega_{n}$ is the measure of $\S^{n-1}$, $\gamma_{\alpha}$ is defined in (\ref{eq:gamma}) and $\overline\gamma_{\alpha}$ in (\ref{eq:gbar})  (compare with (\ref{eq:Lp}) and (\ref{eq:Delta2})). 
In particular,
$S_q^{\rm rad}(\R^{n};\alpha)  =  \omega_{n}^{(q-2)/q}\mu_{q}(\alpha)$, where 
$$
\mu_{q}(\alpha)  =  \inf_{\scriptstyle w\in C^2_c({\R})\atop\scriptstyle w\ne 0}\frac{\displaystyle\int_{-\infty}^\infty\left(|w''|^{2}+2\overline{\gamma}_{\alpha}|w'|^{2}+\gamma_{\alpha}^{2}|w|^{2}\right)ds}{\left(\displaystyle\int_{-\infty}^\infty|w|^q ds\right)^{2/q}}.
$$
If $\gamma_\alpha  =  0$ then clearly $\mu_{q}(\alpha)  =  0$, via rescaling. Conversely, notice that
$\alpha\notin\{4-n,n\}$ if and only if $\gamma_{\alpha}\ne 0$ and in this case the space $H^{2}(\R)$ admits as an equivalent norm
\begin{equation}
\label{eq:w-norm}
\|w\|^{2}_{\alpha}  =  
\int_{-\infty}^\infty\left(|w''|^{2}+2\overline{\gamma}_{\alpha}|w'|^{2}+\gamma_{\alpha}^{2}|w|^{2}\right)ds.
\end{equation}
Thus $\mu_{q}(\alpha)>0$ since $H^2(\R)\hookrightarrow L^q(\R)$ by Sobolev embedding theorem, and hence $S_{q}^{\rm rad}(\R^n;\alpha)>0$.
\QED

We conclude this section with an existence result.

\begin{Theorem}
\label{T:CKNradialexistence}
Let  $q>2$ be a given exponent, and assume that $\alpha\not\in\{4-n,n\}$. 
Then $S_q^{\rm rad}(\R^{n};\alpha)$ is
achieved in $\mathcal{N}^{2}_{\rm rad}(\R^n;\alpha)$.
\end{Theorem}

\proof
Since $\gamma_\alpha\neq 0$ then the Emden-Fowler transform induces an isometry
between $\mathcal N^2_{\rm rad}(\R^n;\alpha)$ and the Sobolev space $H^2(\R)$,
endowed with the equivalent norm in (\ref{eq:w-norm}). It is definitely standard to show the existence of some $\underline{w}\in H^{2}(\R)$ such that $\|\underline{w}\|_{L^{q}}  =  1$ and $\|\underline{w}\|^{2}_{\alpha}  =  \mu_{q}(\alpha)$ (see \cite{Str}). Then the corresponding function $\underline{u}$ defined by (\ref{eq:uwrad}) belongs
to $\mathcal N^2_{\rm rad}(\R^n;\alpha)$ and achieves $S_{q}^{\rm rad}(\R^n;\alpha)$.
\QED

\subsection{Estimates on $S_{q}^{\mathrm{rad}}(\R^{n};\alpha)$}
\label{Ss:esti-Srad}
In this subsection we provide some estimates on the infima
$S_{q}^{\mathrm{rad}}(\R^{n};\alpha)$. We start with the
limiting case $n\ge 5$ and $q  =  \starstar$. 

\begin{Theorem}\label{T:critico-RN}
If $n\ge 5$ and 
\begin{equation}
\label{eq:critico-RN}
2<|\alpha-2|<\sqrt{4+2~\!\frac{(n-2)^{2}(n-4)}{n-3}}
\end{equation}
then $S_{2^{*\!*}}^{\mathrm{rad}}(\R^{n};\alpha)<S^{*\!*}$.
\end{Theorem}

\proof
Set $a  =  -{\alpha}/{2}$. Let $U\in\mathcal{D}^{2}(\R^{n})$ be the radial mapping defined by
\begin{equation}
\label{eq:talenti}
U(x)  =  \left(1+|x|^2\right)^{\frac{4-n}{2}}~\!.
\end{equation}
Our aim is to test $S_{q}^{\mathrm{rad}}(\R^{n};-2a)$
with $|x|^a U$. ~\!In order to simplify notations we put
$$
~\!J  =  \int|x|^{-2}|\nabla U|^{2}~,\qquad ~\!I  =  \int|x|^{-4}|U|^{2}~\!.
$$
We compute
\begin{equation*}
\begin{split}
\int|x|^{-2a}|\Delta(|x|^{a}U)|^{2}&  =  \int|\Delta U|^{2}+4a^{2}~\!J+
a^{2}(n-2+a)^{2}~\!I\\
&\quad+4a^{2}(n-2+a)\int|x|^{-4}U(x\cdot\nabla U)
+4a\int|x|^{-2}(x\cdot\nabla U)\Delta U\\
&\quad +2a(n-2+a)\int|x|^{-2}U\Delta U.
\end{split}
\end{equation*}
Since $U$ is radial, then
\begin{gather*}
\int|x|^{-4}U(x\cdot\nabla U)  =  -\frac{n-4}{2}~\!I~,\quad 
\int|x|^{-2}(x\cdot\nabla U)\Delta U  =  \frac{n}{2}~\!J\\
\int|x|^{-2}U\Delta U  =  -~\!J-(n-4)~\!I.
\end{gather*}
We infer that
$$
\int|x|^{-2a}|\Delta(|x|^{a}U)|^{2}  =  \int|\Delta U|^{2}+2a(a+2)~\!J
-a(a+2)(n-2+a)(n-4-a)~\!I.
$$
Integrating by parts twice we get
$$
J  =  \frac{(n-2)(n-4)^{2}}{4(n-3)}~\!I,
$$
that leads to
$$
\int|x|^{-2a}|\Delta(|x|^{a}U)|^{2}  =  \int|\Delta U|^{2}+a(a+2)\left[a^{2}+2a-\frac{(n-2)^{2}(n-4)}{2(n-3)}\right]~\!I.
$$
Since $U$ achieves the best Sobolev constant $S^{*\!*}$ (see \cite{Sw}), we have
\begin{equation}
\label{eq:US**}
\int|\Delta U|^{2}  =  S^{*\!*}\left(\int |U|^{2^{*\!*}}\right)^{2/2^{*\!*}},
\end{equation}
and then
$$
S_{2^{*\!*}}^{\mathrm{rad}}(\R^{n};\alpha)\le S^{*\!*}+Ca(a+2)\left[a^{2}+2a-\frac{(n-2)^{2}(n-4)}{2(n-3)}\right]I,
$$
where $C>0$ is a power of  the $L^\starstar$ norm of $U$.
Since 
$$
a(a+2)\left[a^{2}+2a-\frac{(n-2)^{2}(n-4)}{2(n-3)}\right]<0
$$
if and only if (\ref{eq:critico-RN}) holds, the conclusion follows.
\QED

Our next goal is to provide the asymptotic behavior of 
$S_{q}^{\mathrm{rad}}(\R^{n};\alpha)$ as $|\alpha|\to \infty$. We first
point out a useful lemma.

\begin{Lemma}
\label{L:equality}
Let $q\ge 2$ and $\alpha,\widetilde\alpha\in\R\setminus\{4-n\}$ be given. Then
$$
S_q^{\rm rad}(\R^n;{\alpha})  =   |\tau(\alpha,\tilde{\alpha})|^{3+\frac{2}{q}}~\inf_{\scriptstyle {  u}\in C^{2}_{c}(\R^n\setminus\{0\})\atop\scriptstyle  u  =  {  u}(|x|),~\scriptstyle  u\ne 0}\frac{\displaystyle\int_{\R^n}|x|^{{\widetilde \alpha}}|\Delta { u}|^2-g(\alpha,{{\widetilde \alpha}})\int_{\R^n}|x|^{\widetilde\alpha-2}|\nabla { u}|^2}
{\left(\displaystyle\int_{\R^n} |x|^{-\beta_{{\widetilde \alpha}}}
|{  u}|^q\right)^{2/q}}~\!,
$$
where
\begin{equation}
\label{eq:tau-g}
\tau(\alpha,\tilde{\alpha}):  =  \dfrac{n-4+ {\alpha}}{n-4+{\widetilde \alpha}}~,\quad g({\alpha},{\widetilde \alpha}):  =  (n-2)~\!\dfrac{({\widetilde\alpha}-{\alpha})[{\widetilde\alpha} {\alpha}-2({\widetilde\alpha}+ {\alpha})-n(n-4)]}{(n-4+{\alpha})^2}.
\end{equation}
\end{Lemma}

\proof
Fix ${u}\in C^\infty_c(\R^n\setminus\{0\})$ radially symmetric, put $\tau  =  \tau(\alpha,\tilde{\alpha})$ and define
\begin{equation}
\label{eq:utau}
{\widetilde u}(r)  =  {u}(r^{1/\tau})~.
\end{equation}
Direct computation leads to
\begin{eqnarray}
\label{eq:taug1}
\int_{\R^n} |x|^{-{\beta_{\alpha}}}|{u}|^q&  =  &|\tau |^{-1}\int_{\R^n} |x|^{-\beta_{\widetilde\alpha}}|{\widetilde u}|^q\\
\nonumber
\int_{\R^n}|x|^{{\alpha}}|\Delta{u}|^2&  =  &
|\tau |^3\int_{\R^n}|x|^{\widetilde\alpha}\left|\Delta {\widetilde u}-(1-\tau ^{-1})(n-2)|x|^{-2}x\cdot\nabla \widetilde u\right|^2 dx\\
\label{eq:taug2}
&  =  &|\tau |^3~\left[\int_{\R^n}|x|^{\widetilde \alpha}|\Delta {\widetilde u}|^2-g({\alpha},{\widetilde \alpha})\int_{\R^n}|x|^{\widetilde \alpha-2}|\nabla {\widetilde u}|^2\right].
\end{eqnarray}
The conclusion readily follows.
\QED

Notice that $g\equiv 0$ in the two-dimensional case. Therefore the
following immediate corollary holds.

\begin{Corollary}\label{C:}
Assume $n  =  2$ and fix  $q\ge 2$. Then for any $\alpha\neq 2$ the ratio
$$
\frac{S_q^{\rm rad}(\R^2;{\alpha})}{|{\alpha}-2|^{3+\frac{2}{q}}}
$$
is a constant, independent on $\alpha$.
\end{Corollary}

Next assume $n\ge 3$.
We will say that ${\alpha}, {\widetilde\alpha}$ are {\em conjugate} if
\begin{equation}
\label{eq:defconj}
({\alpha}-2)({\widetilde\alpha}-2)  =  (n-2)^2.
\end{equation}
Notice that ${\alpha}  =  n$ and ${\widetilde\alpha}  =  4-n$ are self-conjugate. If  ${\alpha}, {\widetilde\alpha}\neq 4-n$ are
conjugate then $g({\alpha},{\widetilde\alpha})  =  0$ and
\begin{equation}
\label{eq:t-conjugate}
|\tau({\alpha}, {\widetilde\alpha}) |  =   \left|\frac{n-2}{{\widetilde\alpha}-2}\right|  =  \left|\frac{{\alpha}-2}{n-2}\right|.
\end{equation}
In this case
\begin{equation}
\label{eq:conjugate}
S_q^{\rm rad}(\R^n;{\alpha})  =  |\tau({\alpha}, {\widetilde\alpha}) |^{3+\frac{2}{q}}S_q^{\rm rad}(\R^n;{\widetilde\alpha}).
\end{equation}

\begin{Corollary}
\label{C:radial2}
Assume $n\ge 3$ and fix $q\ge 2$. Then
$$
\lim_{|\alpha|\to\infty}~\!\frac{S_q^{\rm rad}(\R^n;\alpha)}{|\alpha-2|^{3+\frac{2}{q}}}  =  \frac{S_q^{\rm rad}(\R^n;2)}{(n-2)^{3+\frac{2}{q}}}.
$$
\end{Corollary}

\proof
For any $\alpha\neq 4-n$ we let $\widetilde\alpha$ to be its conjugate
exponent. When $|\alpha|\to \infty$, then $\widetilde\alpha\to 2$ by
(\ref{eq:defconj}), and hence $S_q^{\rm rad}(\R^n;\widetilde\alpha)\to
S_q^{\rm rad}(\R^n;2)$ (use the continuity Lemma \ref{L:cont-alpha} in Appendix
\ref{Ss:continuity}). 
Then the conclusion follows from (\ref{eq:t-conjugate})--(\ref{eq:conjugate}).
\QED

\section{Existence}\label{S:existence}

In this Section we prove Theorems \ref{T:achieved1} and \ref{T:achieved1critical}. We will always
assume:
$$
\left\{\textrm{
\begin{tabular}{l}
$\Sigma\subseteq\S^{n-1}$ is of class $C^2$ and
$-\gamma_{\alpha}\notin\Lambda(\Sigma)$\\
$q>2$  and $q\le 2^{*\!*}$ if $n\ge 5$\\
$\beta  =  n-q~\!\frac{n-4+\alpha}{2}$. 
\end{tabular}}
\right.
$$
In particular,  $S_{q}(\C_\Sigma;\alpha)>0$ by Theorem \ref{T:CKN-Navier}. We need the following result.

\begin{Lemma}
\label{L:minimizer}
Let $(u_{h})\subset\mathcal N^2(\C_\Sigma;\alpha)$ be a minimizing sequence for $S_{q}(\C_\Sigma;\alpha)$. If $u_{h}\weak u$ weakly in $\mathcal N^2(\C_\Sigma;\alpha)$ and $u\ne 0$, then $u$ is a minimizer for $S_{q}(\C_\Sigma;\alpha)$ and $u_{h}\to u$ strongly in $\mathcal N^2(\C_\Sigma;\alpha)$.
\end{Lemma}

The proof is definitely standard. One can adapt to our situation a well known argument (see, e.g., \cite{Str}, Chapt.~1, Sect.~4).

\subsection{$\eps$-compactness}\label{SS:eps-comp}

To prove the existence results stated in the introduction we need an $\eps$-compactness
criterion for sequences of approximating solutions to (\ref{eq:eq}).
We start by pointing out an immediate consequence of Rellich Theorem.

\begin{Lemma}
\label{L:compact}
Let $A$ be a domain with compact closure in $\R^n\setminus\{0\}$. Then $\mathcal N^2(\C_\Sigma;\alpha)$ is compactly embedded into $H^1(\C_\Sigma\cap A)$.
\end{Lemma}


In the next result we let $\mathcal N^{-2}(\C_\Sigma;\alpha)$ to be the topological dual space of $\mathcal N^{2}(\C_\Sigma;\alpha)$ and we use Theorem \ref{T:CKN-Navier} to fix a small number $\eps_0>0$ such that
\begin{equation}
\label{eq:eps}
\eps_0^{\frac{q-2}{q}}<S_q(\C_\Sigma;\alpha)~\!.
\end{equation}

\begin{Proposition}
\label{P:PS}
Let $u_h\in \mathcal N^{2}(\C_\Sigma;\alpha)$, $f_h\in \mathcal N^{-2}(\C_\Sigma;\alpha)$ be
given sequences, such that
$f_h\to 0$ in $\mathcal N^{-2}(\C_\Sigma;\alpha)$, $u_h \weak 0$
 weakly in $\mathcal N^{2}(\C_\Sigma;\alpha)$ and moreover
\begin{gather}
\label{eq:Euler-h}
\Delta\left(|x|^\alpha\Delta u_h\right)  =  |x|^{-\beta}|u_h|^{q-2}u_h+f_h\\
\label{eq:uh-small}
\int_{\C_\Sigma\cap B_R}|x|^{-\beta}|u_h|^q~dx\le \eps_0
\end{gather}
for some $R>0$, where $\eps_0>0$ satisfies (\ref{eq:eps}). Then 
$$
\int_{\C_\Sigma\cap B_{R'}}|x|^{-\beta}|u_h|^q~dx\to 0
\quad\textit{for any $R'\in(0,R)$.}
$$
\end{Proposition}

\proof
Fix $R'\in(0,R)$ and take a cut-off function $\f\in C^\infty_c(B_R)$ such that
$\f\equiv 1$ on $B_{R'}$. Notice that
$$
\int_{\C_\Sigma}|x|^\alpha\Delta u_h\Delta(\f^2 u_h)~dx
  =  \int_{\C_\Sigma}|x|^{\alpha}|\Delta(\f u_h)|^2~dx+o(1)
$$
by Lemma \ref{L:compact}, as $\f$ and its derivatives have compact supports in
$\R^n\setminus\{0\}$. Therefore, from (\ref{eq:Euler-h}) and using
H\"older inequality we get
\begin{eqnarray*}
\int_{\C_\Sigma}|x|^{\alpha}|\Delta(\f u_h)|^2~dx&  =  &
\int_{\C_\Sigma}|x|^{-\beta}|u_h|^{q-2}|(\f u_h)|^2~dx\\
&\le&\left(\int_{\C_\Sigma}|x|^{-\beta}|u_h|^{q}~dx\right)^{\frac{q-2}{q}}~
\left(\int_{\C_\Sigma}|x|^{-\beta}|\f u_h|^{q}~dx\right)^{\frac{2}{q}}.
\end{eqnarray*}
The left hand side in the above inequality can be bounded from below by using the
definition of $S_q(\C_\Sigma;\alpha)$. Thus, from (\ref{eq:uh-small}) we infer
$$
\left(\int_{\C_\Sigma}|x|^{-\beta}|\f u_h|^{q}~dx\right)^{\frac{2}{q}}
S_q(\C_\Sigma;\alpha) \le \eps_0^{\frac{q-2}{q}} 
\left(\int_{\C_\Sigma}|x|^{-\beta}|\f u_h|^{q}~dx\right)^{\frac{2}{q}}~.
$$
The conclusion readily follows from (\ref{eq:eps}), since $\f\equiv 1$ on $B_{R'}$.
\QED

\subsection{Proof of Theorem \ref{T:achieved1}}\label{SS:subcritical}
Using Ekeland's variational principle (see \cite{Str} Chapt.~1, Sect.~5) we can find a minimizing sequence $u_h\in \mathcal N^{2}(\C_\Sigma;\alpha)$, such that (\ref{eq:Euler-h}) holds for a sequence $f_h\to 0$ in $\mathcal N^{-2}(\C_\Sigma;\alpha)$ and such that
$$
\int_{\C_\Sigma}|x|^{\alpha}|\Delta u_h|^2~dx  =  \int_{\C_\Sigma}|x|^{-\beta}|u_h|^q~dx+o(1)
  =   S_q(\C_\Sigma;\alpha)^{\frac{q}{q-2}}+o(1).
$$
Since $u_h$ is bounded in $\mathcal N^{-2}(\C_\Sigma;\alpha)$, we
can assume that $u_h\weak u$ weakly in $\mathcal N^{2}(\C_\Sigma;\alpha)$.
Up to a rescaling, we can also assume that
\begin{equation}
\label{eq:unmezzo}
\int_{\C_\Sigma\cap B_{2}}|x|^{-\beta}|u_h|^q~dx  =  \frac{1}{2} S_q(\C_\Sigma;\alpha)^{\frac{q}{q-2}}.
\end{equation}
We claim that $u\neq 0$. Indeed, if $u_h\weak 0$, then
$$
\int_{\C_\Sigma\cap B_{1}}|x|^{-\beta}|u_h|^q~dx  =  o(1)
$$
by Proposition \ref{P:PS}. On the other hand,
$$
\int_{\C_\Sigma\cap\{1<|x|<2\}}|x|^{-\beta}|u_h|^q~dx  =  o(1)
$$
by Lemma \ref{L:compact} and by Rellich Theorem, contradicting (\ref{eq:unmezzo}).
Thus the minimizing sequence $u_h$ converges weakly to a non trivial limit. Then we can apply Lemma \ref{L:minimizer} to conclude.
\QED

\subsection{Proof of Theorem \ref{T:achieved1critical}}\label{SS:critical}
We put here $S(\alpha)  =  S_{2^{*\!*}}(\C_\Sigma;\alpha)$ to simplify notations. We
select a minimizing sequence $u_h$ as in the proof of Theorem \ref{T:achieved1}. In particular
there exists a sequence $f_h\to 0$ in $\mathcal N^{-2}(\C_\Sigma;\alpha)$ such that
$u_h$ satisfies
\begin{gather}
\label{eq:Euler-h-critical}
\Delta\left(|x|^\alpha\Delta u_h\right)  =  |x|^{\frac{n\alpha}{n-4}}|u_h|^{2^{*\!*}-2}u_h+f_h\\
\label{eq:minimizing}
\int_{\C_\Sigma}|x|^{\alpha}|\Delta u_h|^2~dx  =  \int_{\C_\Sigma}|x|^{\frac{n\alpha}{n-4}}|u_h|^q~dx+o(1)
  =   S(\alpha)^{\frac{n}{4}}+o(1)\\
\label{eq:unmezzo-critical}
\int_{\C_\Sigma\cap B_{2}}|x|^{\frac{n\alpha}{n-4}}|u_h|^{2^{*\!*}}~dx
  =  \frac{1}{2} S(\alpha)^{\frac{n}{4}}.
\end{gather}
As before, we have to prove that $u_h$ cannot converge weakly to $0$. By contradiction,
assume that $u_h\weak 0$ weakly in $\mathcal N^{2}(\C_\Sigma;\alpha)$.
Then we can argue as in the proof of Theorem \ref {T:achieved1} to get
$$
\int_{\C_\Sigma\cap B_{1}}|x|^{\frac{n\alpha}{n-4}}|u_h|^{2^{*\!*}}~dx  =  o(1)
$$
and hence
\begin{equation}
\label{eq:contradiction-critical}
\int_{\C_\Sigma\cap\{1<|x|<2\}}|x|^{\frac{n\alpha}{n-4}}|u_h|^{2^{*\!*}}~dx  =  \frac{1}{2} S(\alpha)^{\frac{n}{4}}
+o(1)
\end{equation}
by (\ref{eq:unmezzo-critical}). Now we take a cut-off function $\f\in C^\infty_c(\R^n\setminus\{0\})$
such that $\f\equiv 1$ on $B_2\setminus B_1$ and we use $\f^2 u_h$ as test function in (\ref{eq:Euler-h-critical}).
Using Lemma \ref{L:compact}, H\"older inequality and (\ref{eq:minimizing}) we have
$$
\int_{\C_\Sigma}|x|^{\alpha}|\Delta(\f u_h)|^2~dx
\le
S(\alpha)~
\left(\int_{\C_\Sigma}|x|^{\frac{n\alpha}{n-4}}|\f u_h|^{2^{*\!*}}~dx\right)^{\frac{n-4}{n}}+o(1).
$$
Let us define
$$
F_h  =  \Delta(|x|^{\alpha/2}\f u_h)-|x|^{\alpha/2}\Delta(\f u_h).
$$
Using Lemma \ref{L:compact} and Rellich Theorem one plainly gets that
$F_h\to 0$ strongly in $L^2(\R^n)$. Thus, by Sobolev inequality,
\begin{eqnarray*}
\int_{\C_\Sigma}|x|^{\alpha}|\Delta(\f u_h)|^2~dx&  =  &
\int_{\C_\Sigma}|\Delta(|x|^{\frac{\alpha}{2}}\f u_h)|^2~dx+o(1)\\
&\ge& S^{*\!*}\left(\int_{\C_\Sigma}\left||x|^{\frac{\alpha}{2}}\f u_h\right|^{2^{*\!*}}~dx\right)^{\frac{n-4}{n}}+o(1).
\end{eqnarray*}
Putting together these informations we conclude that
$$
S^{*\!*}\left(\int_{\C_\Sigma}|x|^{\frac{n\alpha}{n-4}}|\f u_h|^{2^{*\!*}}~dx\right)^{\frac{n-4}{n}}
\le S(\alpha)~
\left(\int_{\C_\Sigma}|x|^{\frac{n\alpha}{n-4}}|\f u_h|^{2^{*\!*}}~dx\right)^{\frac{n-4}{n}}+o(1).
$$
Thus
$$
o(1)  =  \int_{\C_\Sigma}|x|^{\frac{n\alpha}{n-4}}|\f u_h|^{2^{*\!*}}~dx
\ge \int_{\C_\Sigma\cap\{1<|x|<2\}}|x|^{\frac{n\alpha}{n-4}}|\f u_h|^{2^{*\!*}}~dx,
$$
as $0<S(\alpha)<S^{*\!*}$ by assumption and $\f\equiv 1$ on the annulus $B_2\setminus B_1$. Since this conclusion contradicts (\ref{eq:contradiction-critical}), we infer that the weak limit
of the minimizing sequence $u_h$ cannot vanish. Then we can apply Lemma \ref{L:minimizer} to conclude.
\QED 

{From} Theorems \ref{T:achieved1critical}, \ref{T:minore} and  \ref{T:critico-RN} we infer the next existence result.

\begin{Theorem}
\label{T:achieved1critical2}
Let $n\ge 5$ and let $\Sigma$ be a domain in $\S^{n-1}$ of class $C^2$. Assume that $-\gamma_{\alpha}\notin\Lambda(\Sigma)$. The best constant $S_{2^{*\!*}}(\mathcal C_\Sigma;\alpha)$ is achieved if one of the following conditions holds:
\begin{itemize}
\item[$(i)$] $n\ge 6$ and $|\alpha-2|>2$
\item[$(ii)$] $n = 5$, $\Sigma=\S^{4}$ and $2<|\alpha-2|<\sqrt{13}$.
\end{itemize}
\end{Theorem}

\begin{Remark}
Assume $\Sigma = \S^{n-1}$ and $n\ge 6$.
By Proposition \ref{P:large} it results that $S_{2^{*\!*}}(\R^n\setminus\{0\};\alpha)\le
S^{*\!*}$, while $S_\starstar^{\rm rad}(\R^n;\alpha)$ diverges as $|\alpha|\to\infty$.
Thus, for $|\alpha|$ large enough, extremals for $S_{2^{*\!*}}(\R^n\setminus\{0\};\alpha)$ 
do exist, 
but none of them is radially symmetric. This breaking symmetry phenomenon is definitively 
new with respect to the Caffarelli-Kohn-Nirenberg first order inequalities.
Breaking symmetry will be studied in more detail in Section
\ref{S:BS}.
\end{Remark}

\section{Breaking positivity}\label{S:BP}

In this section we illustrate a surprising phenomenon that is completely new
with respect to similar first order problems. Namely, we show that all functions achieving
the best constant $S_{q}(\mathcal{C}_{\Sigma};\alpha)$ might be forced to change sign. 
In particular,
extremal functions for $S_{q}(\mathcal{C}_{\Sigma};\alpha)$ cannot be positive if $q$ is close to 2 and 
\begin{equation}
\label{eq:break-pos}
-\gamma_{\alpha}>\frac{\lambda_{1}+\lambda_{2}}{2}~\!,
\end{equation}
where $\lambda_{1}$ and $\lambda_{2}$ are the two first eigenvalues of $-\Delta_{\sigma}$ in $H^{1}_{0}(\Sigma)$. To this goal we introduce the infima
\begin{equation*}
S_{q}^{+}(\mathcal{C}_{\Sigma};\alpha)  =  \inf_{\scriptstyle u\in C^{2}(\overline{\mathcal{C}_{\Sigma}}\setminus\{0\})\atop\scriptstyle u\ge 0,~u\not\equiv 0}\frac{\displaystyle\int_{\mathcal{C}_{\Sigma}}|x|^{\alpha}|\Delta u|^{2}dx}{\displaystyle\left(\int_{\mathcal{C}_{\Sigma}}|x|^{-\beta}u^{q}dx\right)^{2/q}}.
\end{equation*}
Let us state the main result of this section.

\begin{Theorem}
\label{T:break-pos}
Assume 
(\ref{eq:break-pos}). Then there exists $q_{\alpha}>2$ such that 
$$S_{q}(\mathcal{C}_{\Sigma};\alpha)<S_{q}^{+}(\mathcal{C}_{\Sigma};\alpha)$$ 
for all $q\in[2,q_{\alpha})$. In particular, if $q\in(2,q_{\alpha})$, extremal functions for $S_{q}(\mathcal{C}_{\Sigma};\alpha)$ cannot be positive.
\end{Theorem}

\proof
In order to prove Theorem \ref{T:break-pos} we will use once more the Emden-Fowler
transform $T$ already introduced in (\ref{eq:EF}). Besides the infimum
$m_{N}(\Sigma;\gamma)$ in (\ref{eq:mN}), we define also
$$
m_{N}^{+}(\Sigma;\gamma)  =  \inf_{\scriptstyle\varphi\in H^{2}\cap H^{1}_{0}(\Sigma)\atop\scriptstyle\varphi\ge 0,~\varphi\not\equiv 0}\frac{\displaystyle\int_{\Sigma}|-\Delta_{\sigma}\varphi+\gamma\varphi|^{2}d\sigma}{\displaystyle\int_{\Sigma}\varphi^{2}d\sigma}~.
$$
The following facts hold:
\begin{itemize}
\item[$(i)$]
$S^+_2(\mathcal C_\Sigma;\alpha) = m^+_N(\Sigma;\gamma_\alpha)$.
\item[$(ii)$] 
If $-\gamma\le\frac{\lambda_{1}+\lambda_{2}}{2}$ then $m_{N}^{+}(\Sigma;\gamma)  =  m_{N}(\Sigma;\gamma)$.
\item[$(iii)$] 
If $-\gamma>\frac{\lambda_{1}+\lambda_{2}}{2}$ then $m_{N}^{+}(\Sigma;\gamma)>m_{N}(\Sigma;\gamma)$. In particular $m_{N}^{+}(\Sigma;\gamma)>0$ for all $\gamma\in\R$, $\gamma\ne\lambda_{1}$.
\end{itemize}
Theorem \ref{T:break-pos} is an immediate consequence of $(i)-(iii)$ and of the continuity Lemma
\ref{L:cont-q} in Appendix \ref{Ss:continuity}.
Claim $(i)$ easily follows from the computations in Section \ref{S:Inequalities} on the Emden-Fowler transform.
To prove
$(ii)$, notice that if $-\gamma\le\frac{\lambda_{1}+\lambda_{2}}{2}$ then $\mathrm{dist}(-\gamma,\Lambda(\Sigma))  =  |\gamma+\lambda_{1}|$. Thus, by $(i)$ in Lemma \ref{L:aux}, $m_{N}(\Sigma;\gamma)$ is achieved by an eigenfunction $\varphi_{1}$ corresponding to $\lambda_{1}$. Since one can take $\varphi_{1}\ge 0$, it follows that $m_{N}^{+}(\Sigma;\gamma)\le m_{N}(\Sigma;\gamma)$. The opposite inequality is trivial.\\
Now we check $(iii)$. Assume that $m_{N}^{+}(\Sigma;\gamma)\le m_{N}(\Sigma;\gamma)$. Then equality holds. By a standard argument one can plainly check that $m_{N}^{+}(\Sigma;\gamma)$ is attained, namely there exists $\varphi\in  H^{2}\cap H^{1}_{0}(\Sigma)$ such that 
$$
\int_{\Sigma}|-\Delta_{\sigma}\varphi+\gamma\varphi|^{2}d\sigma  =  m_{N}^{+}(\Sigma;\gamma)~\!,\quad\int_{\Sigma}\varphi^{2}d\sigma  =  1~\!,\quad\varphi\ge 0~\!.
$$
Then $\varphi$ is an extremal for $m_{N}(\Sigma;\gamma)$, too. By $(i)$, $\varphi$ is an eigenfunction of $-\Delta_{\sigma}$. Since $\varphi\ge 0$, it must be $-\Delta_{\sigma}\varphi  =  \lambda_{1}\varphi$ and, again by $(i)$, $\mathrm{dist}(-\gamma,\Lambda(\Sigma))  =  |\gamma+\lambda_{1}|$, that is, $-\gamma\le\frac{\lambda_{1}+\lambda_{2}}{2}$.
Theorem \ref{T:break-pos} is completely proved.
\QED

Specializing Theorem \ref{T:break-pos} to the case $\Sigma  =  \S^{n-1}$, when $\lambda_{1}  =  0$ and $\lambda_{2}  =  n-1$, we immediately obtain the next result.

\begin{Corollary}
Assume that 
$$
|\alpha-2|>\sqrt{(n-1)^{2}+1}.
$$
Then there exists $q_{\alpha}>2$ such that $S_{q}(\R^{n}\setminus\{0\};\alpha)<S_{q}^{+}(\R^{n}\setminus\{0\};\alpha)$ for all $q\in[2,q_{\alpha})$. In particular, if $q\in(2,q_{\alpha})$, extremal functions for $S_{q}(\R^{n}\setminus\{0\};\alpha)$ cannot be positive.
\end{Corollary}

\section{Breaking symmetry}\label{S:BS}

In this section we discuss some conditions for breaking symmetry. We use the the constants $S_q(\R^n\setminus\{0\};\alpha)$ and $S_q^{\rm rad}(\R^n;\alpha)$ already defined in (\ref{eq:SqRnalpha}) and (\ref{eq:minimization1radial}), respectively.

As a first condition, we have that if $-\gamma_{\alpha}$ is close enough to the spectrum $\Lambda(\S^{n-1})$ then breaking symmetry occurs.

\begin{Theorem}
\label{T:break-symm1}
For every $q>2$ and for every $k\in\mathbb{N}$ there exists $\delta>0$ such that if 
$$
0<|\gamma_{\alpha}+k(n-2+k)|<\delta
$$ 
then $S_q(\R^n\setminus\{0\};\alpha)<S_q^{\rm rad}(\R^n;\alpha)$.
\end{Theorem}

\proof
Fix $k\in\mathbb{N}$, let $\lambda  =  k(n-2+k)$ and $\alpha_{0}$ be such that $-\gamma_{\alpha_{0}}  =  \lambda$. Since $\lambda\in\Lambda(\S^{n-1})$, by Theorem \ref{T:CKN-Navier} it turns out that
$S_q(\R^n\setminus\{0\};\alpha_{0})  =  0$. In general, if $\alpha\to\alpha_{0}$ then 
$$
S_q(\R^n\setminus\{0\};\alpha_{0})\ge\limsup_{\alpha\to\alpha_{0}}S_q(\R^n\setminus\{0\};\alpha)
$$
(see Remark \ref{R:usc}). Hence we have that
$S_q(\R^n\setminus\{0\};\alpha)\to 0\quad\textrm{as }\alpha\to\alpha_{0}$. We also have 
$S_q^{\mathrm{rad}}(\R^n;\alpha)\to S_q^{\mathrm{rad}}(\R^n;\alpha_{0})$ as $\alpha\to\alpha_{0}$
by Lemma \ref{L:cont-alpha},
and $S_q^{\mathrm{rad}}(\R^n;\alpha_{0})>0$ by Theorem \ref{T:CKNradial}. Hence the conclusion follows.
\QED

As a second condition, we show that if $|\alpha|$ is large then again breaking symmetry occurs. More precisely we have the following result.


\begin{Theorem}
\label{T:bs}
Let $q>2$, $q\le 2^{*\!*}$ when $n\ge 5$, and let $\alpha\in\R$.
If 
\begin{equation}
\label{eq:bs1}
|\gamma_{\alpha}|>\frac{n-1}{q-2}\left(1+\sqrt{q-1}\right)
\end{equation}
then $S_q(\R^n\setminus\{0\};\alpha)<S_q^{\rm rad}(\R^n;\alpha)$.
\end{Theorem}

\proof
Assume that $S_q(\R^n\setminus\{0\};\alpha)=S_q^{\rm rad}(\R^n;\alpha)$ for some
$\alpha\in\R$ such that $\gamma_\alpha\neq 0$. We claim that in this case
\begin{equation}
\label{eq:bs}
|\gamma_\alpha|\le\frac{n-1}{q-2}~\!\left(1+\sqrt{q-1}\right)~\!.
\end{equation}
We start by noticing that 
$S_q^{\rm rad}(\R^n;\alpha)>0$  by Theorem \ref{T:CKNradial}. 
Thus also $S_q(\R^n\setminus\{0\};\alpha)>0$, and hence $-\gamma_{\alpha}$ does not belong to the spectrum of the Laplace-Beltrami operator on $\S^{n-1}$, by Theorem \ref{T:CKN-Navier}. In particular, the Hilbert space $\mathcal N^2(\R^n\setminus\{0\};\alpha)$ is well defined. 
We introduce the following functionals on $\mathcal N^2(\R^n\setminus\{0\};\alpha)\setminus\{0\}$:
$$
A(u):  =  \int|x|^\alpha|\Delta u|^2~,\quad B(u):  =  \left(\int|x|^{-\beta}|u|^q\right)^{2/q}~,
\quad
R(u):  =  \frac{A(u)}{B(u)}~\!.
$$
Let $\uu$ be the radially symmetric solution to the minimization problem
(\ref{eq:minimization1radial}) given by Theorem \ref{T:CKNradial}. Thus $\uu$
achieves also $S_q(\R^n\setminus\{0\};\alpha)$, that is, $\uu$ minimizes the functional $R(u)$ on 
$\mathcal N^2(\R^n\setminus\{0\};\alpha)\setminus\{0\}$. In particular,
\begin{equation}
\label{eq:R}
R'(\uu)\cdot v  =  0~,\quad R''(\uu)[v,v]\ge 0\quad
\textrm{for any $v\in \mathcal N^2(\R^n\setminus\{0\};\alpha)$.}
\end{equation}
In order to simplify notation we can assume that $B(\uu)  =  1$. Then by
direct computations based on (\ref{eq:R}) one gets
\begin{gather}
\label{eq:minimum}
B''(\uu)[v,v]~\!\left(\int|x|^\alpha|\Delta \uu|^2\right)\le A''(\uu)[v,v]  =  2\int|x|^\alpha|\Delta v|^2\\
\nonumber
B''(\uu)[v,v]  =  2\left[(q-2)\left(\int|x|^{-\beta}|\uu|^{q-1}\uu v\right)^2+
(q-1)\int|x|^{-\beta}|\uu|^{q-2}|v|^2~\right].
\end{gather}
Now we choose the test function $v$, that is, $v  =  \underline{u}\f$, where $\f\in H^1(\S^{n-1})$
is an eigenfunction of the Laplace-Beltrami operator on the sphere, relatively
to the first positive eigenvalue and normalized with respect to the  $L^2$ norm. Hence
$$
-\Delta_\sigma\f  =  (n-1)\f~,\quad \int_{\S^{n-1}}\f  =  0~,\quad\int_{\S^{n-1}}|\f|^2  =  1.
$$
Since $\f$ has zero mean value, then
$$
B''(\uu)[\uu\f,\uu\f]  =  2
(q-1)\int|x|^{-\beta}|\uu|^{q}  =  2(q-1).
$$
Then we compute
\begin{eqnarray*}
\int|x|^\alpha|\Delta (\uu\f)|^2&  =  &
\int|x|^\alpha|(\Delta \uu -(n-1)|x|^{-2}\uu)|^2\\
&  =  &
\int|x|^\alpha|\Delta \uu|^2+(n-1)^2\int|x|^{\alpha-4}|\underline{u}|^2
-2(n-1)\int|x|^{\alpha-2}\uu\Delta\uu\\
&\le&
\int|x|^\alpha|\Delta \uu|^2+(n-1)^2\int|x|^{\alpha-4}|\underline{u}|^2\\
&&\quad\quad\quad
+2(n-1)\left(\int|x|^{\alpha-4}|\uu|^2\right)^{1/2}
\left(\int|x|^\alpha|\Delta\underline{u}|^2\right)^{1/2}
\end{eqnarray*}
by the Cauchy-Schwarz inequality. Thus from (\ref{eq:minimum}) we get
\begin{eqnarray*}
(q-2)\!\int|x|^\alpha|\Delta \uu|^2&\le& (n-1)^2\int|x|^{\alpha-4}|\uu|^2\\
&&\quad\quad\quad+
2(n-1)\left(\int|x|^{\alpha-4}|\uu|^2\right)^{\!1/2}\!\!
\left(\int|x|^\alpha|\Delta\uu|^2\right)^{\!1/2}.
\end{eqnarray*}
Thus
$(q-2)\xi^2\le(n-1)^2+2(n-1)\xi$, where
$$
\xi :=  \left(\frac{\displaystyle\int|x|^\alpha|\Delta\uu|^2}{\displaystyle
\int|x|^{\alpha-4}|\uu|^2}\right)^{1/2}\ge ~\!\ |\gamma_\alpha|
$$
by (\ref{eq:inf-radial}).
Inequality (\ref{eq:bs}) readily follows
via elementary calculus.
\QED

\begin{Remark}
Assume $n\le 4$ and fix any $\alpha\notin\{4-n,n\}$. If $q>2$ is
large enough then  breaking symmetry occurs,
that is, $S_q(\R^n\setminus\{0\};\alpha)<S_q^{\mathrm{rad}}(\R^n;\alpha)$.
\end{Remark}

\section{Dirichlet boundary conditions}\label{S:Dbc}

In this section we assume that $\Sigma$ is a domain of class $C^2$ with compact
closure in $\S^{n-1}$. In particular, 
$\partial \mathcal C_\Sigma\setminus\{0\}$ is not empty. Our aim is to study 
study minimization problems of the form (\ref{eq:dirichlet}).  
First of all we recall that
\begin{equation}\label{eq:lin-dir}
\inf_{\scriptstyle u\in C^{2}_{c}(\mathcal C_\Sigma)\atop\scriptstyle u\ne 0}\frac{\displaystyle\int_{\mathcal C_\Sigma}|x|^{\alpha}|\Delta u|^{2}dx}{\displaystyle\int_{\mathcal C_\Sigma}|x|^{\alpha-4}|u|^{2}dx}>0~\!,
\end{equation}
see \cite{CM1}, whatever $\alpha\in\R$ is. Thus we can define the Hilbert space $\mathcal{D}^{2}(\mathcal{C}_{\Sigma};\alpha)$ as the completion of $C^{2}_{c}(\mathcal{C}_{\Sigma})$ with respect to the norm defined in (\ref{eq:norm}). 
In particular, if $q>2$, $q\le 2^{*\!*}$ if $n\ge 5$, by density we have that 
\begin{equation}
\label{eq:SqD}
S_q^D(\mathcal C_\Sigma;\alpha)=\inf_{\scriptstyle u\in \mathcal{D}^{2}(\mathcal C_\Sigma;\alpha)\atop\scriptstyle u\ne 0}\frac{\displaystyle\int_{\mathcal C_\Sigma}|x|^{\alpha}|\Delta u|^{2}dx}{\displaystyle\left(\int_{\mathcal C_\Sigma}|x|^{-\beta}|u|^{q}dx\right)^{2/q}}~\!.
\end{equation}

Our main results about the existence of minimizers for problems (\ref{eq:SqD}) are summarized in the next theorem.

\begin{Theorem}\label{T:Dir}
Let $\alpha\in\R$ and let $\Sigma$ be a domain of class $C^2$ properly contained in $\S^{n-1}$.
Let $q>2$, and $q\le \starstar$ if $n\ge 5$. Then $S_q^{D}(\mathcal C_\Sigma;\alpha)>0$ and
moreover:
\begin{itemize}
\item[$(i)$]
If $n\le 4$ or $q<2^{*\!*}$ then $S_q^{D}(\mathcal C_\Sigma;\alpha)$ is achieved in $\mathcal{D}^{2}(\mathcal{C}_{\Sigma};\alpha)$.
\item[$(ii)$]
For $n\ge 5$, if $S_{2^{*\!*}}^{D}(\mathcal C_\Sigma;\alpha)<S^{*\!*}$, then $S_{2^{*\!*}}^{D}(\mathcal C_\Sigma;\alpha)$ is achieved in $\mathcal{D}^{2}(\mathcal{C}_{\Sigma};\alpha)$.
\item[$(iii)$]
If $n\ge 6$, then $S_{2^{*\!*}}^{D}(\mathcal C_\Sigma;\alpha)<S^{*\!*}$.
\item[$(iv)$]
If $S_{q}^{D}(\mathcal C_\Sigma;\alpha)$ is attained, then $S_{q}^{N}(\mathcal C_\Sigma;\alpha)<S_{q}^{D}(\mathcal C_\Sigma;\alpha)$.
\end{itemize}
\end{Theorem}

\proof
Parts $(i)$--$(iii)$ can be proved by repeating the same argument developed in Sections \ref{S:Inequalities} and \ref{S:existence}.
As far as concerns $(iv)$, we point out that the large inequality $S_{q}^{N}(\mathcal C_\Sigma;\alpha)\le S_{q}^{D}(\mathcal C_\Sigma;\alpha)$ always holds true. Moreover, if $u\in\mathcal{D}^{2}(\mathcal C_\Sigma;\alpha)$ is a minimizer for $S_{q}^{D}(\mathcal C_\Sigma;\alpha)$ and equality holds, $u$  would be a solution of 
$$
\Delta(|x|^{\alpha}\Delta u)  =  \lambda|x|^{-\beta}|u|^{q-2}u~~\textrm{in }\mathcal{C}_{\Sigma}
$$
for some $\lambda>0$ and it would satisfy both Neumann and Dirichlet boundary conditions. Hence it would be $u  =  0$, which is impossible.\QED

\appendix
\section{Remarks on a Brezis-Nirenberg type problem}\label{S:BN}

In this section we deal with the Dirichlet problem
\begin{equation}
\label{eq:BN-problemA}
\begin{cases}
\Delta^2 u+\lambda\Delta u  =  |u|^{\starstar-2}u&\textrm{in $B$}\\
u  =  u(|x|)~,~u\neq 0\\
u  =  |\nabla u|  =  0&\textrm{on $\partial B$}
\end{cases}
\end{equation}
where $\lambda\in\R$ is a given parameter and $B$ is the unit ball in $\R^n$, $n\ge 5$.
Since a detailed analysis of problem (\ref{eq:BN-problemA}) would lead us far from our
purposes, we limit ourself to investigate those features of problem (\ref{eq:BN-problemA}) 
that have some relevance with the questions under investigation in the present paper.

We point out that the fourth order differential
equation in (\ref{eq:BN-problemA}) contains a leading term with critical growth and a linear term involving the Laplacian. In the spirit of the result by Brezis-Nirenberg \cite{BreNir}, this last term provides a perturbation of a dilation invariant problem which allows us to recover compactness, when the parameter $\lambda$ stays in a suitably restricted range. 

We start our analysis by pointing out a non-existence result.

\begin{Theorem}\label{T:BN5}
If $\lambda\le 0$ then problem (\ref{eq:BN-problemA}) has no solution.
If $n  =  5$ and $\lambda\le 21/8$ then problem (\ref{eq:BN-problemA}) has no solution.
\end{Theorem}

\proof
For $\lambda  =  0$ the result is already known, see for instance 
\cite{Gaz98} or \cite{GGS}. If $\lambda\neq 0$,
the proof is based on a Pohozaev identity that has to be coupled with a Hardy-type inequality
in the lowest
dimensional case.

Let $u$ be a solution of (\ref{eq:BN-problemA}).
We put $r  =  |x|$ and we denote by $u_{r}$ the radial derivatives of $u$, namely $u_{r}  =  r^{-1}x\cdot\nabla u$. Testing (\ref{eq:BN-problemA}) with $2ru_{r}-u$ one infers the following
Pohozaev identity (use for instance the computations in \cite{GGS}, pagg. 250--252):
\begin{equation}
\label{eq:A}
2\lambda\int_B|\nabla u|^2  =  \omega_n\int_{\partial B}|\Delta u|^2(x\cdot \nu)  =  \omega_n |u_{rr}(1)|^2,
\end{equation}
where $\omega_n$ is the measure of $\S^{n-1}$. Thus $\lambda>0$.

Now we assume $n  =  5$. We will prove in a moment that
\begin{equation}
\label{eq:r3}
5\int_B r^2|\Delta u|^2-6\int_B|\nabla u|^2-2\lambda\int_B|\nabla u|^2  =  -\lambda \int_B r^2|\nabla u|^2-
\frac{7}{5}\int_B r^2|u|^{10}.
\end{equation}
{From} (\ref{eq:r3}) and using Lemma \ref{L:AAA} with $\alpha = 2$ we get
$$
0> 5\int_B r^2|\Delta u|^2-6\int_B|\nabla u|^2-2\lambda\int_B|\nabla u|^2
\ge\left(\frac{21}{4}-2\lambda\right)\int_B|\nabla u|^2,
$$
that implies $\lambda>21/8$ and concludes the proof.
\\
It remains to check (\ref{eq:r3}). If $\eta$ and $\f$ are radial and smooth enough, then
\begin{equation}
\label{eq:eta}
\int_B \f\eta_{r}  =  \omega_4\f(1)\eta(1)-\int_B(\f_{r}+4r^{-1}\f)\eta~,\quad
2\int_B\eta u_{r}\Delta u  =  \int_B(4r^{-1}\eta-\eta_{r})|\nabla u|^2.
\end{equation}
Next we notice that
$$
{u}\Delta {u}  =  \frac{1}{2}\Delta({u}^2)-|\nabla {u}|^2~,\quad u_{rr}  =  \Delta {u}-4r^{-1}u_{r}~,\quad
(\Delta {u})_{r}  =  \Delta u_{r}-4r^{-2}u_{r},
$$
and we test (\ref{eq:BN-problemA})
with $r^3 u_{r}$. Using integration by parts, (\ref{eq:eta})  and (\ref{eq:A}) we get
\begin{equation*}\begin{split}
\int_B(\Delta^2 u)~\!(r^3 u_{r})& =  -\omega_4|u_{rr}(1)|^2+
\int_B(\Delta u)(r^3\Delta u_{r}+u_{r}\Delta(r^3)+6r^2u_{rr})\\
&  =  -2\lambda\int_B|\nabla u|^2+\int_B r^3
(\Delta u)(\Delta u)_{r}+6\int_B r^2|\Delta u|^2
-18 \int_B ru_{r}\Delta u\\
&  =  -\lambda\int_B|\nabla u|^2+\frac{5}{2} \int_B r^2|\Delta u|^2
-3\int_B |\nabla u|^2.
\end{split}\end{equation*}
We compute also
$$
\lambda\int_B\Delta u(r^3 u_{r})  =  \frac{\lambda}{2}\int_B r^2|\nabla u|^2~\!,\quad
\int_B u^9(r^3 u_{r})  =  \frac{1}{{10}}\int_B r^3\left(|u|^{10}\right)_{r}  =  -\frac{7}{{10}}\int_B r^2 |u|^{10}~\!.
$$
Thus, from (\ref{eq:BN-problemA}) we readily get (\ref{eq:r3}).
\QED

A natural approach for studying (\ref{eq:BN-problemA}) consists
in looking for minimizers for
\begin{equation}
\label{eq:BN}
S^{*\!*}_{\lambda}:  =  
\inf_{\scriptstyle u\in H^{2}_{0,\mathrm{rad}}(B)\atop\scriptstyle u\ne 0}
\frac{\displaystyle\int_{B}|\Delta u|^{2}dx-\lambda\int_{B}|\nabla u|^{2}dx}{\displaystyle\left(\int_{B}|u|^\starstar dx\right)^{2/2^{*\!*}}}~\!.
\end{equation}
Clearly, the infimum $S^{*\!*}_{\lambda}$ is positive provided that $\lambda<\lambda_{2,1}$,
where 
$$
\lambda_{2,1}:  =  \inf_{\scriptstyle u\in H^{2}_{0,\mathrm{rad}}(B)\atop\scriptstyle u\ne 0}\frac{\displaystyle\int_{B}|\Delta u|^{2}dx}{\displaystyle\int_{B}|\nabla u|^{2}dx}
\ge\frac{n^2}{4}
$$
by Lemma \ref{L:AAA} in Appendix \ref{Ss:continuity}. Moreover, minimizers for
$S^{*\!*}_{\lambda}$ give rise to solutions to problem (\ref{eq:BN-problemA}).
Arguing for instance as in Proposition \ref{P:large} one can check that
$S^{*\!*}_{\lambda}\le S^{*\!*}$
for any $\lambda\in \R$. In particular, by  monotonicity, it turns out
that $S_{\lambda}^{*\!*}  =  S^{*\!*}$ and is not attained if $\lambda\le 0$, accordingly
with Theorem \ref{T:BN5}. As in 
\cite{BreNir} or \cite{PLL84}, a crucial point in finding an existence result for
the minimization problem (\ref{eq:BN}) consists in giving sufficient conditions
for the validity of the strict inequality $S_{\lambda}^{*\!*}<S^{*\!*}$. First of all
we notice that $S_{\lambda}^{*\!*}<S^{*\!*}$ provided that $\lambda$  is close enough to $\lambda_{2,1}$. Notice indeed that for any $\lambda>0$ it results
$$
S^{*\!*}_{\lambda}\le
\inf_{\scriptstyle u\in H^{2}_{0,\mathrm{rad}}(B)\atop\scriptstyle u\ge 0~\!, ~\!u\ne 0}\frac{\displaystyle\int_{B}|\Delta u|^{2}dx-{\lambda}{\lambda_1^{-1}}\int_{B}| u|^{2}dx}{\displaystyle\left(\int_{B}|u|^{\starstar}dx\right)^{2/2^{*\!*}}}~\!,
$$
where $\lambda_1$ is the Poincar\'e constant of the unit ball in $\R^n$. Then one concludes
by using known results for problem
$$
\begin{cases}
\Delta^2 u-\lambda u  =  |u|^{\starstar-2}u&\textrm{in $B$}\\
u  =  |\nabla u|  =  0&\textrm{on $\partial B$}
\end{cases}
$$
(see for instance \cite{Gaz98} or \cite{GGS}). The same argument shows
that $S^{*\!*}_{\lambda}<S^{*\!*}$ if $n\ge 8$. The next lemma, that was crucially used in Section \ref{SS:Limiting}, covers also the case $n\in\{6,7\}$ and shows that $n = 5$ is the only
critical dimension for problem (\ref{eq:BN-problemA}).

\begin{Lemma}
\label{L:BN-esti}
Let $B$ be the unit ball in $\R^{n}$ and $\lambda>0$. If $n\ge 6$, then there exists
a nonnegative radially symmetric function $u\in C^\infty_c(B)$ such that
$$
\frac{\displaystyle\int_{B}|\Delta u|^{2}dx-{\lambda}
\int_{B}|\nabla u|^{2}dx}{\displaystyle\left(\int_{B}|u|^{\starstar}dx\right)^{2/2^{*\!*}}}<S^{*\!*}.
$$
\end{Lemma}

\proof
Let $U$ be the non-negative radial mapping defined in (\ref{eq:talenti}) and let
$\xi\in C^{\infty}_{c}(B)$ be a radial function with $0\le\xi\le 1$ and $\xi(x)  =  1$ as $|x|\le\frac{1}{2}$. Define
$$
u_{\eps}(x)  =  \eps^{\frac{4-n}{2}}\xi(x)U(\eps^{-1}x).
$$
Hence $u_{\eps}\in C^{2}_{c}(B)$ and 
\begin{equation}
\label{eq:ueps1}
\int|\Delta u_{\eps}|^{2}  =  \int|\Delta U|^{2}+O(\eps^{n-4})\quad\textrm{and}\quad\int u_{\eps}^{2^{*\!*}}  =  \int U^{2^{*\!*}}+O(\eps^{n})\quad\textrm{as $\eps\to 0$}
\end{equation}
(see, e.g., \cite{Gaz98}).
Thanks to (\ref{eq:US**}) and (\ref{eq:ueps1}) we have that
\begin{equation}
\label{eq:ueps2}
\int|\Delta u_{\eps}|^{2}  =  \left(S^{*\!*}+O(\eps^{n-4})\right)\left(\int U^{2^{*\!*}}\right)^{2/2^{*\!*}}\quad\textrm{as $\eps\to 0$}.
\end{equation}
If $n\ge 7$ then $U\in D^{1,2}(\R^{n})$ and one can easily check that
\begin{equation*}
\int|\nabla u_{\eps}|^{2}  =  \eps^{2}\int|\nabla U|^{2}+o(\eps^{2})\quad\textrm{as $\eps\to 0$.}
\end{equation*}
Therefore 
\begin{equation*}
\frac{S^{*\!*}_{\lambda}-S^{*\!*}}{\eps^{2}}\le
\frac{1}{\eps^{2}}\left[\frac{\displaystyle\int|\Delta u_{\eps}|^{2}}{\displaystyle\left(\int u_{\eps}^{2^{*\!*}}\right)^{2/2^{*\!*}}}-S^{*\!*}\right]-\lambda\frac{\displaystyle\int|\nabla U|^{2}}{\displaystyle\left(\int U^{2^{*\!*}}\right)^{2/2^{*\!*}}}+o(1)\quad\textrm{as $\eps\to 0$}
\end{equation*}
and then $S^{*\!*}_{\lambda}<S^{*\!*}$, because of (\ref{eq:ueps2}). 
If $n  =  6$ then
$$
\frac{1}{\eps^{2}}\int|\nabla u_{\eps}|^{2}\ge\int_{B_{\frac{1}{2\eps}}}|\nabla U|^{2}  =  C\int_{0}^{\frac{1}{2\eps}}\frac{r^{n-1}}{(1+r^{2})^{n-2}}~\!dr\ge C|\log\eps|
$$
for some constant $C>0$. Then 
$$
\frac{S^{*\!*}_{\lambda}-S^{*\!*}}{\eps^{2}}\le
\frac{1}{\eps^{2}}\left[\frac{\displaystyle\int|\Delta u_{\eps}|^{2}}{\displaystyle\left(\int u_{\eps}^{2^{*\!*}}\right)^{2/2^{*\!*}}}-S^{*\!*}\right]-\lambda\frac{\displaystyle{\eps^{-2}}\int|\nabla u_{\eps}|^{2}}{\displaystyle\left(\int u_{\eps}^{2^{*\!*}}\right)^{2/2^{*\!*}}}\le O(1)-\lambda\frac{C|\log\eps|}{\displaystyle\left(\int u_{\eps}^{2^{*\!*}}\right)^{2/2^{*\!*}}}.
$$
Hence also in this case we can conclude that $S^{*\!*}_{\lambda}<S^{*\!*}$.
\QED

We conclude this section with an existence result, whose proof can be obtained
by using the above remarks and standard arguments.

\begin{Theorem}\label{T:BN6}
Let $B$ be the unit ball in $\R^{n}$ and $0<\lambda<\lambda_{2,1}$. 
\begin{itemize}
\item[$(i)$] If $n\ge 6$ then 
problem (\ref{eq:BN-problemA}) admits a ground state solution, i.e., a function $u\in H^{2}_{0}(B)$ solving (\ref{eq:BN-problemA}) and minimizing $S^{*\!*}_{\lambda}$.
\item[$(ii)$] If $n  =  5$ then there exists $\lambda^*\in(0,\lambda_{2,1})$ such that
$S^{*\!*}_{\lambda}$ is not achieved if  $\lambda<\lambda^*$ and achieved if $\lambda^*<\lambda<\lambda_{2,1}$. 
\end{itemize}
\end{Theorem}

\section{Auxiliary results and open problems}
\label{S:Appendix}

This Appendix contains some technical results used in the previous sections. In particular we 
prove of some estimates that were used in the proof of Theorem
\ref{T:minore} and a couple of continuity lemmas. Finally we write a list of open problems.

\begin{Lemma}
\label{L:radial-esti}
Let $a\in\R$ and $e\in\mathbb{S}^{n-1}$. Then there exists a constant $K_{a}>0$ such that for every radial mapping $u\in C^{2}_{c}(B)$ and for every $t\in[0,1]$ one has
$$
\int|tx+e|^{-2a}\left|\Delta\left(|tx+e|^{a}u\right)\right|^{2}\le
\int|\Delta u|^{2}-2C_{a}t^{2}\int|\nabla u|^{2}+K_{a}t^{3}\int|\nabla u|^{2}
$$
where $C_{a}  = {a(a+2)(n-2)}/{n}$.
\end{Lemma}
%


\proof
One computes
$$
\Delta\left(|tx+e|^{a}u\right)  =  |tx+e|^{a-2}\left[|tx+e|^{2}\Delta u+2at\nabla u\cdot(tx+e)+a(n-2+a)t^{2}u\right]
$$
and then
\begin{equation*}
\begin{split}
&\int|tx+e|^{-2a}\left|\Delta\left(|tx+e|^{a}u\right)\right|^{2}  =  
\int|\Delta u|^{2}+
4at\underbrace{\int|tx+e|^{-2}(\nabla u\cdot(tx+e))\Delta u}_{I_{1}}\\
&~+4a^{2}t^{2}\underbrace{\int|tx+e|^{-4}|\nabla u\cdot(tx+e)|^{2}}_{I_{2}}+2a(n-2+a)t^{2}\underbrace{\int|tx+e|^{-2}u\Delta u}_{I_{3}}\\
&~+4a^{2}(n-2+a)t^{3}\!\underbrace{\int|tx+e|^{-4}u(\nabla u\cdot(tx+e))}_{I_{4}}+a^{2}(n-2+a)^{2}t^{4}\!\underbrace{\int|tx+e|^{-4}|u|^{2}}_{I_{5}}.
\end{split}
\end{equation*}
Since $u$ is radial one has that
\begin{equation}
\label{eq:radial-case}
\begin{array}{c}
\displaystyle\int|\nabla u\cdot e|^{2}  =  \int\frac{(e\cdot x)^{2}}{|x|^{2}}~\!|\nabla u|^{2}  =  \frac{1}{n}\int|\nabla u|^{2}\vspace{6pt}\\
\displaystyle\int\frac{e\cdot x}{|x|^{2}}~\!|\nabla u|^{2}  =  \int(x\cdot\nabla u)(e\cdot\nabla u)  =  0.
\end{array}
\end{equation}
For future convenience we also point out that since
\begin{gather}
\nonumber
\nabla |tx+e|^{-2}  =  -2t|tx+e|^{-4}(e+x)\\
\label{eq:A3}
1-|tx+e|^{-2}  =  t|tx+e|^{-2}(2x\cdot e+t|x|^{2})
\end{gather}
and since $|tx+e|\ge 1-t$ for $x\in B$, we have the following estimates:
\begin{equation}
\label{eq:IIIII}
|tx+e|^{-4}\le C~,\quad
\left|\nabla |tx+e|^{-2}\right|\le Ct~,\quad
\left|1-|tx+e|^{-2}\right|\le Ct
\end{equation}
for all $x\in B$ and for every $t\ge 0$ small enough.
\medskip

\noindent
\textbf{Estimate of $I_{1}$.}
Firstly we integrate by parts, obtaining that
\begin{equation*}
\begin{split}
I_{1}&  =  \int|tx+e|^{-2}(\nabla u\cdot(tx+e))\Delta u  =  -\int\nabla\left(|tx+e|^{-2}(tx+e)\cdot\nabla u\right)\cdot\nabla u\\
&  =  2t^{2}\int|tx+e|^{-4}\left[(x\cdot\nabla u)^{2}+(e\cdot\nabla u)(x\cdot\nabla u)\right]\\
&\quad+2t\int|tx+e|^{-4}\left[(e\cdot\nabla u)^{2}+(e\cdot\nabla u)(x\cdot\nabla u)\right]\\
&\quad+\frac{n-2}{2}t\int|tx+e|^{-2}|\nabla u|^{2}-t^{4}\int|tx+e|^{-4}\left(e\cdot x+|x|^{2}\right)|\nabla u|^{2}\\
&\quad-\int|tx+e|^{-2}\left[\frac{e\cdot x}{|x|}u_{rr}u_{r}+\frac{e\cdot\nabla u}{|x|}u_{r}-\frac{e\cdot x}{|x|^{2}}u_{r}^{2}\right]
\end{split}
\end{equation*}
where, as in the proof of Theorem \ref{T:BN5}, $u_{r}$ and $u_{rr}$ denote the first and second radial derivatives of $u$, respectively.
Then we use the fact that $u$ is radial, in particular the identity $e\cdot\nabla u  =  \frac{e\cdot x}{|x|}~\!u_{r}$, and (\ref{eq:radial-case}), getting that
\begin{equation*}
\begin{split}
I_{1}&  =  \frac{n-1}{2}\underbrace{\int|tx+e|^{-2}~\!\frac{e\cdot x}{|x|^{2}}~\!|\nabla u|^{2}}_{J_{1}}
+\left(\frac{n}{2}-1\right)t\underbrace{\int|tx+e|^{-2}|\nabla u|^{2}}_{J_{2}}\\
&\quad+t\underbrace{\int|tx+e|^{-4}\left[e\cdot x+\frac{(e\cdot x)^{2}}{|x|^{2}}
\right]|\nabla u|^{2}}_{J_{3}}
+t^{2}\underbrace{\int|tx+e|^{-4}\left[e\cdot x+|x|^{2}\right]|\nabla u|^{2}}_{J_{4}}.
\end{split}
\end{equation*}
We need to compute the terms of $I_{1}$ of order 0 and 1 in $t$. Therefore, in view of (\ref{eq:radial-case}) and (\ref{eq:A3}), we can write
\begin{equation*}
\begin{split}
J_{1}&  =  \int\left[|tx+e|^{-2}-1\right]\frac{e\cdot x}{|x|^{2}}~\!|\nabla u|^{2}\\
&  =  -t^{2}\int|tx+e|^{-2}(e\cdot x)|\nabla u|^{2}-2t\int|tx+e|^{-2}~\!\frac{(e\cdot x)^{2}}{|x|^{2}}~\!|\nabla u|^{2}.
\end{split}
\end{equation*}
Next, using again (\ref{eq:radial-case}) and (\ref{eq:A3}), we estimate
\begin{equation*}
\begin{split}
\int|tx+e|^{-2}&~\!\frac{(e\cdot x)^{2}}{|x|^{2}}~\!|\nabla u|^{2}  =  \int|tx+e|^{-2}\left[\frac{(e\cdot x)^{2}}{|x|^{2}}-1\right]|\nabla u|^{2}+\int\frac{(e\cdot x)^{2}}{|x|^{2}}~\!|\nabla u|^{2}\\
&  =  -t\int|tx+e|^{-2}\left[2e\cdot x+t|x|^{2}\right]\frac{(e\cdot x)^{2}}{|x|^{2}}~\!|\nabla u|^{2}+\frac{1}{n}\int|\nabla u|^{2}.
\end{split}
\end{equation*}
Hence, by (\ref{eq:IIIII}),
$$
J_{1}\le-\frac{2t}{n}\int|\nabla u|^{2}+Ct^{2}\int|\nabla u|^{2}.
$$
Then, using (\ref{eq:IIIII}) too, we estimate
$$
J_{2}  =  \int|\nabla u|^{2}+\int\left[|tx+e|^{-2}-1\right]|\nabla u|^{2}\le\int|\nabla u|^{2}+Ct\int|\nabla u|^{2}.
$$
Moreover, again by (\ref{eq:radial-case}), (\ref{eq:A3}) and (\ref{eq:IIIII}), we have 
\begin{equation*}
\begin{split}
J_{3}&  =  
\int\left[|tx+e|^{-4}-1\right]\left[e\cdot x+\frac{(e\cdot x)^{2}}{|x|^{2}}\right]|\nabla u|^{2}+\int\left[e\cdot x+\frac{(e\cdot x)^{2}}{|x|^{2}}\right]|\nabla u|^{2}\\
&  =  -t\int|tx+e|^{-2}\left[|tx+e|^{-2}+1\right]\left[e\cdot x+\frac{(e\cdot x)^{2}}{|x|^{2}}\right]\left[2e\cdot x+t|x|^{2}\right]|\nabla u|^{2}\\
&\quad+\frac{1}{n}\int|\nabla u|^{2}\le\frac{1}{n}\int|\nabla u|^{2}+Ct\int|\nabla u|^{2}
\end{split}
\end{equation*}
and
$$
J_{4}\le C\int|\nabla u|^{2}.
$$
In conclusion
$$
I_{1}\le\left(\frac{2}{n}+\frac{n}{2}-2\right)t\int|\nabla u|^{2}+Ct^{2}\int|\nabla u|^{2}
$$
for $t\ge 0$ small enough.
\medskip

\noindent
\textbf{Estimate of $I_{2}$.}
Thanks to (\ref{eq:radial-case}) and (\ref{eq:IIIII}), we have
\begin{equation*}
\begin{split}
I_{2}&  =  \int|tx+e|^{-4}|\nabla u\cdot(tx+e)|^{2}\\
&\le(1-t)^{-4}\left[t^{2}\int|x\cdot\nabla u|^{2}+2t\int(x\cdot\nabla u)(e\cdot\nabla u)+\int|e\cdot\nabla u|^{2}\right]\\
&\le Ct\int|\nabla u|^{2}+\frac{1}{n}\int|\nabla u|^{2}
\end{split}
\end{equation*}
for $t\ge 0$ small enough.
\medskip

\noindent
\textbf{Estimate of $I_{3}$.}
Firstly we integrate by parts, then we use Cauchy-Schwarz and Poincar\'e inequalities and the estimates (\ref{eq:IIIII}) getting that
\begin{equation*}
\begin{split}
I_{3}&  =  \int|tx+e|^{-2}u\Delta u  =  -\int\nabla\left(|tx+e|^{-2}u\right)\cdot\nabla u\\
&  =  2t\int|tx+e|^{-4}u~\!(x+e)\cdot\nabla u+\int\left(1-|tx+e|^{-2}\right)|\nabla u|^{2}-\int|\nabla u|^{2}\\
&\le 4t(1-t)^{-4}\int|u|~\!|\nabla u|+Ct\int|\nabla u|^{2}-\int|\nabla u|^{2}\\
&\le Ct\int|\nabla u|^{2}-\int|\nabla u|^{2}
\end{split}
\end{equation*}
for $t\ge 0$ small enough.
\medskip

\noindent
\textbf{Estimate of $I_{4}$.}
Using Cauchy-Schwarz and Poincar\'e inequalities we have that
$$
I_{4}  =  \int|tx+e|^{-4}u(\nabla u\cdot(tx+e))\le(1-t)^{-4}\int|u|~\!|\nabla u|\le C\int|\nabla u|^{2}
$$
for $t\ge 0$ small enough.
\medskip

\noindent
\textbf{Estimate of $I_{5}$.} Thanks to the Poincar\'e inequality we have
$$
I_{5}  =  \int|tx+e|^{-4}|u|^{2}\le(1-t)^{-4}\int|u|^{2}\le C\int|\nabla u|^{2}
$$
for $t\ge 0$ small enough.
\medskip

\noindent
The conclusion easily follows from the previous estimates of $I_{1}$,...,$I_{5}$.\QED

\subsection{Continuity lemmas}\label{Ss:continuity}

Here we discuss the following two auxiliary results.

\begin{Lemma}\label{L:cont-q}
If $q\to 2^{+}$ then $S_{q}(\mathcal{C}_{\Sigma};\alpha)\to S_{2}(\mathcal{C}_{\Sigma};\alpha)$ and $S_{q}^{+}(\mathcal{C}_{\Sigma};\alpha)\to S_{2}^{+}(\mathcal{C}_{\Sigma};\alpha)$.
\end{Lemma}

\begin{Lemma}\label{L:cont-alpha}
If $\alpha\to \alpha_{0}$ then $S_{q}^{\mathrm{rad}}(\R^{n};\alpha)\to S_{q}^{\mathrm{rad}}(\R^{n};\alpha_{0})$.
\end{Lemma}

Let us point out the following general fact.

\begin{Remark}\label{R:usc}
Let $X$ be a nonempty set, let $I$ be an interval in $\R$ and let $F\colon X\times I\to [0,+\infty)$ be a given mapping. For every $a\in I$ set $S(a)  =  \inf_{u\in X}F(u,a)$. If $F$ is continuous with respect to the parameter $a\in I$, then $S(a)\ge\limsup S(a_{k})$ when $a_{k}\to a$. This fact can be proved in an elementary way. 
\end{Remark}

In view of Remark \ref{R:usc}, in order to prove Lemmas \ref{L:cont-q} and \ref{L:cont-alpha} we need to show just the lower semicontinuity inequalities. 

\subsubsection*{Proof of Lemma \ref{L:cont-q}}
It is a consequence of Remark \ref{R:usc} and of the following general result.

\begin{Lemma}
Let $\Omega$ be a domain in $\R^{n}$ and let $\mu$ be a positive measure on $\Omega$. Let $X$ be a space of measurable functions from $\Omega$ into $\R$, endowed with some norm $\|\cdot\|$. Assume that $X$ is continuously embedded into $L^{p}(\mu)$ for all $p$ in some compact interval $I\subset[1,\infty)$. Then the mapping
$$
p\mapsto S(p):  =  \inf_{\scriptstyle{u\in X}\atop\scriptstyle{u\ne 0}}\frac{\|u\|}{|u|_{p}}\quad\textrm{where }|u|_{p}  =  \left(\int_{\Omega}|u|^{p}~\!d\mu\right)^{\frac{1}{p}}\textrm{ and }p\in  I
$$
is lower semicontinuous in $I$, i.e., if $(p_{k})\subset I$ and $p_{k}\to p$ then $S(p)\le\liminf S(p_{k})$.
\end{Lemma}

\proof
Let $I  =  [p_{0},p_{1}]$. By H\"older's inequality, for every $u\in X$ and $\theta\in[0,1]$ one has that
$$
|u|_{p_{\theta}}^{p_{\theta}}\le|u|_{p_{0}}^{(1-\theta)p_{0}}|u|_{p_{1}}^{\theta p_{1}}
$$
where $p_{\theta}  =  \theta p_{1}+(1-\theta)p_{0}$. This readily implies that
\begin{equation}
\label{eq:S(q)conc}
S(p_{\theta})^{p_{\theta}}\ge S(p_{0})^{(1-\theta)p_{0}}S(p_{1})^{\theta p_{1}}.
\end{equation}
Setting $f(p)  =  p\log S(p)$, (\ref{eq:S(q)conc}) reads
$$
(1-\theta) f(p_{0})+\theta f(p_{1})\le f(\theta p_{1}+(1-\theta)p_{0})
$$
namely $f$ is concave in $[p_{0},p_{1}]$. This implies that $f$ is continuous in $(p_{0},p_{1})$ and lower semicontinuous in $[p_{0},p_{1}]$. Clearly the same holds for $S(p)$, too. \QED

\subsubsection*{Proof of Lemma \ref{L:cont-alpha}}

It is trivial if $n=2$, by Corollary \ref{C:}. Thus assume $n\ge 3$. 
If $\alpha_{0}\in\{n,4-n\}$ then $S_{q}^{\mathrm{rad}}(\R^{n};\alpha_{0})  =  0$ (see Theorem \ref{T:CKNradial}) and the result is a consequence of Remark \ref{R:usc}. If $\alpha_{0}\in\R\setminus\{n,4-n\}$ then the proof of Lemma \ref{L:cont-alpha} can be accomplished according to the following argument. First of all let us introduce the infimum
$$
\mu_{2,1}^{\mathrm{rad}}(\R^{n};\alpha):  =  \inf_{\scriptstyle u\in C^{2}_{c}(\R^{n}\setminus\{0\})\atop\scriptstyle u  =  u(|x|),~u\ne 0}\frac{\displaystyle\int_{\R^{n}}|x|^{\alpha}|\Delta u|^{2}~\!dx}{\displaystyle\int_{\R^{n}}|x|^{\alpha-2}|\nabla u|^{2}~\!dx}~\!.
$$

\begin{Lemma}\label{L:AAA}
There results $\displaystyle\mu_{2,1}^{\mathrm{rad}}(\R^{n};\alpha)  =  \left(\frac{n-\alpha}{2}\right)^{2}.$
\end{Lemma}

\proof
Proceeding as in the proof of Theorem \ref{T:CKNradial}, by means of the Emden-Fowler transform we have that
$$
\mu_{2,1}^{\mathrm{rad}}(\R^{n};\alpha)-\left(\frac{n-\alpha}{2}\right)^{2}  =  \inf_{\scriptstyle w\in C^2_c({\R})\atop\scriptstyle w\ne 0}
\frac{\displaystyle\int_{-\infty}^\infty|w''|^{2}ds+\left[2\overline{\gamma}_{\alpha}-\left(\frac{n-\alpha}{2}\right)^{2}\right]\int_{-\infty}^\infty|w'|^{2}ds}
{\displaystyle\int_{-\infty}^\infty|w'|^{2}ds+\left(\frac{n-4+\alpha}{2}\right)^{2}\int_{-\infty}^\infty|w|^{2}ds}
$$
where $\overline{\gamma}_{\alpha}$ is defined in (\ref{eq:gbar}). We notice that 
$$
2\overline{\gamma}_{\alpha}-\left(\frac{n-\alpha}{2}\right)^{2}  =  \left(\frac{n-4+\alpha}{2}\right)^{2}
$$
and we conclude by a standard scaling argument.
\QED

\begin{Lemma}\label{L:alal}
If $\alpha,\tilde{\alpha}\in\R\setminus\{n,4-n\}$ then
\begin{equation}
\label{eq:alal}
\left[1-\frac{4|g(\alpha,\tilde{\alpha})|}{(n-\tilde\alpha)^{2}}\right]S_{q}^{\mathrm{rad}}(\tilde\alpha)\le |\tau(\tilde{\alpha},\alpha)|^{3+\frac{2}{q}}S_{q}^{\mathrm{rad}}(\alpha)\le\left[1+\frac{4|g(\alpha,\tilde{\alpha})|}{(n-\alpha)^{2}}\right]S_{q}^{\mathrm{rad}}(\tilde\alpha)
\end{equation}
where $\tau(\tilde\alpha,\alpha)$ and $g(\alpha,\tilde{\alpha})$ are defined in (\ref{eq:tau-g}).
\end{Lemma}

\proof
As in the proof of Lemma \ref{L:equality}, for every ${u}\in C^\infty_c(\R^n\setminus\{0\})$ radially symmetric let $\widetilde{u}\in C^\infty_c(\R^n\setminus\{0\})$ be the radial function defined by means of the transformation (\ref{eq:utau}). Thanks to the identities (\ref{eq:taug1})--(\ref{eq:taug2}) and using the definition of $\mu_{2,1}^{\mathrm{rad}}(\R^{n};\alpha)$ and Lemma \ref{L:AAA}, the conclusion readily follows.
\QED
\medskip

\noindent
\textit{Completion of the proof of Lemma \ref{L:cont-alpha}.} Let $\alpha_{0}\in\R\setminus\{n,4-n\}$ and $\alpha_{k}\to\alpha_{0}$. We can apply Lemma \ref{L:alal} and, since $\tau(\alpha_{0},\alpha_{k})\to 1$ and $g(\alpha_{k},\alpha_{0})\to 0$, the conclusion follows from (\ref{eq:alal}).
\QED

\subsection{Remarks and open problems}\label{Ss:open}

The arguments in Section \ref{S:existence} can be used to get existence results of
minimizers for
$$
S_{q,\lambda}(\mathcal C_\Sigma;\alpha)  =  \inf_{\scriptstyle u\in \mathcal{N}^{2}_\lambda(\mathcal{C}_{\Sigma};\alpha)\atop\scriptstyle u\ne 0}\frac{\displaystyle\int_{\mathcal C_\Sigma}|x|^{\alpha}|\Delta u|^{2}dx
-\lambda\int_{\mathcal C_\Sigma}|x|^{\alpha-4}|u|^2dx}
{\displaystyle\left(\int_{\mathcal C_\Sigma}|x|^{-\beta}|u|^{q}dx\right)^{2/q}}~\!.
$$
Here $\alpha\in \R$ is any parameter, $q>2$ and $q\le\starstar$ if $n\ge 5$, 
$\lambda<\textrm{dist}\left(-\gamma_{\alpha},\Lambda(\Sigma)\right)^{2}$ (compare with \cite{CM1}),
and $\mathcal{N}^{2}_\lambda(\mathcal{C}_{\Sigma};\alpha)$ is a suitably defined function space.
The approach in Section \ref{S:existence} can be plainly applied also 
to prove existence results for extremals of Lin's inequality in \cite{Lin86} and
for more general dilation-invariant inequalities.

\medskip

The present paper raises several open questions. We list few of them.

\medskip

i) It might be interesting to generalize the results of this paper when $|\Delta u|^{2}$ is replaced by $|\Delta u|^{p}$ with $p>1$. Some partial results can be found in \cite{SzWa}.

\medskip

ii) 
Our results about breaking positivity and breaking symmetry hold only for some restricted ranges of $\alpha$ and/or $q$. Is it possible to give a sharper description of the region of parameters $\alpha, q$ where
breaking positivity/symmetry occur? In particular, is it true that breaking positivity occurs
for any $\alpha$ large enough?

\medskip

iii)  
Is it true that for $n\ge 3$ and $\alpha\in(4-n,n)$ extremals for $S_{q}(\R^{n}\setminus\{0\};\alpha)$ are radially symmetric and/or positive?

\medskip

iv)  
Let $\Sigma$ be properly contained in $\S^{n-1}$, and take $n\ge 5$, $\alpha=0$. Then
$$
0<S_\starstar(\mathcal C_\Sigma;0)\le S_\starstar^D(\mathcal C_\Sigma;0)=S^{*\!*}~\!.
$$
Is it true that $S_\starstar(\mathcal C_\Sigma;0)=S^{*\!*}$? This question
 is related to 
\cite{GazGruSwe}.

\medskip

v) 
When $n\ge 6$ we showed that $S_{2^{*\!*}}(\mathcal C_{\Sigma};\alpha)<S^{*\!*}$ if $|\alpha-2|>2$. Indeed we suspect that if  $\alpha\in(0,4)$ then $S_{2^{*\!*}}(\mathcal C_{\Sigma};\alpha)=S^{*\!*}$ and is not achieved.

\label{References}

\footnotesize

\bigskip\bigskip
\noindent
Paolo Caldiroli\\
Dipartimento di Matematica\\
Universit\`a di Torino, via Carlo Alberto, 10\\
10123 Torino, Italy. \\
Email: {paolo.caldiroli@unito.it}

\bigskip

\noindent
Roberta Musina\\
{Dipartimento di Matematica ed Informatica}\\
Universit\`a di Udine, via delle Scienze, 206\\
33100 Udine, Italy. \\
Email: {roberta.musina@uniud.it}

\end{document}